\documentclass{amsart}
   \usepackage[cp850]{inputenc}
  \RequirePackage{amssymb} % for the square for the proof environment.
  \RequirePackage{amsmath} % for multiple lines in equations
  \usepackage{amsmath}
  \usepackage[mathscr]{eucal}
 \usepackage{color}

  \usepackage{amsfonts}
  \usepackage{latexsym}
   \RequirePackage{amsbsy} % for Boldsymbol
   \RequirePackage{amsopn} % for DeclareMathOperator
   \RequirePackage{amsmath}% for multiple lines in equations
    \RequirePackage{amssymb}% for subsetneq

 \DeclareMathOperator{\co} {\boldsymbol{c}}

    \DeclareMathOperator{\calF}     {\mathcal F}%
  \DeclareMathOperator{\calM}     {\mathcal M}%
  \DeclareMathOperator{\calS}     {\mathcal S}%
  \DeclareMathOperator{\calT}     {\boldsymbol{\mathcal T}}%
   \DeclareMathOperator{\calbS}     {\boldsymbol{\mathcal S}}%
  \DeclareMathOperator{\calV} {\boldsymbol{\mathcal V}}%
  \DeclareMathOperator{\calW}     {\boldsymbol{\mathcal W}}%
    \DeclareMathOperator{\calL} {\boldsymbol{\mathcal L}}%

\DeclareMathOperator{\calbT} {\boldsymbol{\mathcal T}}

   \DeclareMathOperator{\Kr} {\mbox{\rm Kr}}%

\DeclareMathOperator{\Na} {\mbox{\rm Na}}

\DeclareMathOperator{\Spec} {\mbox{\rm Spec}}
\DeclareMathOperator{\Max} {\mbox{\rm Max}}
\DeclareMathOperator{\KN} {\mbox{\rm KN}}
\DeclareMathOperator{\SKN} {\mbox{\rm KN$^\prime$}}

\DeclareMathOperator{\QSpec}
{\mbox{\rm QSpec}} \newtheorem{thee}{Theorem}[section]
\newtheorem{coor}[thee]{Corollary} \newtheorem{leem}[thee]{Lemma}

\newtheorem{prro}[thee]{Proposition}
\newtheorem{deef}[thee]{Definition}
\newtheorem{exxe}[thee]{Example} \newtheorem{reem}[thee]{Remark}
\newtheorem{t-d}[thee]{Theorem/Definition}

   \newcommand{\balf}
   {\renewcommand{\theenumi}{(\alph{enumi})}
   \renewcommand{\labelenumi}{\theenumi}
                        \begin{enumerate}}
  \newcommand{\ealf}   {\end{enumerate}
                        \renewcommand{\theenumi}{\arabic{enumi}}
                        \renewcommand{\labelenumi}{\theenumi.}}
  \newcommand{\bara}   {\renewcommand{\theenumi}{(\arabic{enumi})}
                        \renewcommand{\labelenumi}{\theenumi}
                        \begin{enumerate} }
  \newcommand{\eara}   {\end{enumerate}
                        \renewcommand{\theenumi}{\arabic{enumi}}
                        \renewcommand{\labelenumi}{\theenumi.}}

   \newcommand{\brom} {\renewcommand{\theenumi}{
(\roman{enumi})} \renewcommand{\labelenumi}{\theenumi}
                        \begin{enumerate} }

  \newcommand{\erom}   {\end{enumerate}
                        \renewcommand{\theenumi}{\arabic{enumi}}
                        \renewcommand{\labelenumi}{\theenumi.}}

\begin{document}

\title {A generalization of Kronecker function rings and Nagata rings}

\medskip

 \author{ Marco Fontana}
\address{Dipartimento di Matematica
Universit\`a degli Studi Roma Tre \\
Largo San Leonardo Murialdo, 1 \\
00146 Roma, Italy}
\email{fontana@mat.uniroma3.it}

\author{Alan Loper}
\address{Department of Mathematics \\ Ohio State University-Newark\\ 
Newark,
Ohio 43055 USA}
\email{lopera@math.ohio-state.edu}

\begin{abstract} Let $D$ be an integral domain with quotient field $K$.  The Nagata ring 
$D(X)$ and the Kronecker function ring Kr$(D)$ are both subrings of the field of rational functions $K(X)$ containing as a subring
the ring $D[X]$ of polynomials in the variable $X$.  Both of these function
rings have been extensively studied and generalized.  The principal 
interest in these two extensions of $D$ lies in the reflection of 
various algebraic and spectral properties of $D$  and Spec$(D)$ in algebraic and spectral properties  of the function 
rings.  Despite the obvious similarities in definitions and properties, 
these two kinds of domains of rational functions have been classically treated 
independently, when $D$ is not a Pr\"ufer domain.  The purpose of this note is to study two different 
unified approaches to the Nagata rings and the Kronecker function 
rings, which yield these rings and their classical generalizations as 
special cases.  \end{abstract}


 \maketitle

\medskip 
%\centerline{OCTOBER 14, 2005  }

%%%%%%%%%%%%%%%%%% SECTION 1
\section{Introduction}

Let $D$ be a  commutative   integral domain with
quotient field $K$.  Let $X$ be an indeterminate over $D$ and let 
$f \in D[X]$.  We denote by 
$\boldsymbol{c}(f)$ \it the content of the polynomial $f$, \rm  i.e.  
$\boldsymbol{c}(f)$ is the ideal of $D$
generated by the coefficients of $f$. Moreover, if $\calV(D)$ is the 
set of all the valuation overrings of $D$, for each ideal $I$ of $D$, 
we set  $I^{b} := \bigcap \{IV \mid V \in \calV(D)\}$ (cf. \cite[page 
398]{G} and \cite[Appendix 4]{ZS}).

Two classical rings related to $D$ which have both been well studied
are the Nagata ring $D(X)$ and the Kronecker function ring $\Kr(D)$
defined as follows.
%\begin{deef} \label{def:1.1}
$$\Na(D) := D(X) := \left\{ \frac{f}{g} \mid  f,g \in D[X],\ 
\boldsymbol{c}(f)
\subseteq \boldsymbol{c}(g), \ \boldsymbol{c}(g) \mbox{ \rm  is 
invertible} \right\}
$$
 is \it the Nagata ring of $D$. \rm \
%\end{deef}
(Note that this is not the most common definition:  $D(X)$ is usually
defined by designating that $f,g \in D[X]$ and that 
$\boldsymbol{c}(g) = D$,  \cite[Section 33]{G}.    The
definition above is equivalent to this one and fits our program
better.)

 On the other hand, if $D$ is integrally closed, then:  
$$\Kr(D) := \left\{ \frac{f}{g} \mid f,g \in D[X],\ g\neq 0,\ 
\boldsymbol{c}(f)^{b} \subseteq \boldsymbol{c}(g)^{b} \right\} $$
 is \it the Kronecker function ring of
$D$,  \cite[Section 32]{G}. \rm

These two rings of rational functions are the same if and only if $D$ 
is a
Pr\"ufer domain   \cite[Theorem 33.4]{G}.    In fact, both rings 
arose as generalizations of
Kronecker's original function rings which specified that $D$ should
be  a ring of algebraic numbers or, more generally,   a Dedekind 
domain (and, hence, a Pr\"ufer domain),  (cf. \cite{K}, \cite{W} 
and \cite {E}).  When $D$ is any arbitrary integrally closed 
domain it is easy to see that $\Na(D) \subseteq \Kr(D)$, 
\cite[Theorem 33.3]{G}.

There are obvious similarities in the two definitions in spite of the
fact that they generally yield different rings.  Next we give
equivalent definitions  (or, characterizations)    for each type 
of ring in which there are also
obvious similarities,  i.e. both these rings can be constructed by 
intersection of families of Nagata rings of quasilocal overrings.  \ 
(Note: We do assume for both of these results that we
know how to construct the Nagata ring  $R(X)$  for a \sl 
quasilocal \rm domain $R$; in this situation, the condition 
``$\boldsymbol{c}(g)$  is invertible''  becomes  ``$ 
\boldsymbol{c}(g)$  is principal'', \cite[Proposition 7.4]{G}.)

%THEOREM/DEFINITION  1.1
\begin{t-d} \label{t-d:1} \rm \cite[Theorem 33.3, Theorem 32.10 and 
the proof of Corollary 32.14]{G}  \sl Let $D$ be an integral domain, 
let $\Max(D)$ [respectively, $ \Spec(D)$]  represent the set of all 
maximal [respectively, prime]  ideals of
$D$ and let  $\calV(D)$ [respectively, $\calV_{\mbox{\footnotesize 
\it \!min}}(D)$] denote the set of all the valuation [respectively, 
minimal valuation] overrings of $D$.
\begin{enumerate}

\item[\bf (1)] \sl  $\Na(D) = \bigcap \left\{ D_{M}(X) \mid M \in 
\Max(D) \right\} = \bigcap \left\{ D_{P}(X) \mid P \in \Spec(D) 
\right\}$\,.

\item[\bf (2)] \sl  If $D$ is integrally closed,  $\Kr(D) = \bigcap 
\left\{ V(X) \mid V \in \calV(D) \right\} = \bigcap \left\{ W(X) \mid 
\right.$ $\left.W \in \calV_{\mbox{\footnotesize \it \!min}}(D) 
\right\}$\,.
\end{enumerate}
\end{t-d}

Both the Kronecker function ring and the Nagata ring have been
generalized and intensively
studied   (cf. for instance \cite{Kr2}, \cite{Arnold:1969}, 
\cite{Arnold-Brewer:1971},  \cite{A:1977}, \cite{Kang:1989}, 
\cite{OM3},  \cite{FL1}, \cite{FL3} and\cite{HK3}).    However, in 
spite of their common origin, they have been
studied separately.  There are generalized Nagata rings and
generalized Kronecker function rings which are distinct objects of
study.  A major goal of this paper is to define and study a single
construction for a class of function rings which includes the
Kronecker function ring, the Nagata ring and their generalizations as
special cases.

Following the double characterization of the Kronecker function rings
and Nagata rings above   (Theorem/Definition \ref{t-d:1}),    we 
approach our generalized function rings from
two separate directions, as rings of individually chosen rational
functions and as intersections of Nagata rings of  quasilocal 
    overrings.

In generalizing the rational function approach,  we note that the
standard generalization of Kronecker function rings,   introduced 
by Krull \cite{Kr2},    involves replacing
the $b$--operation  with a more general  ``star--operation''   
(for short, $\star$--operation)   belonging to a
special class of star operations known as the ``e.a.b. star
operations''   (the explicit definitions are recalled in Section 
2).    The point is that e.a.b. operations have some nice
properties in common with the $b$--operation which make possible
the proof that the Kronecker function ring,    as defined by Krull 
(cf. the definition after  Remark \ref{rk:3.4}),    is actually a
ring  for any arbitrary integrally closed domain.

Note that in the definition of the Kronecker
function ring, any  nonzero   polynomial is eligible to be the 
denominator of a
rational function.   Note further that we
do not allow arbitrary  nonzero   polynomials in the 
denominators for the Nagata
ring.  Rather, we only allow a very restricted class of polynomials,
those with invertible content.
 Finally note that   the definition of the Nagata ring
given is  formally   comparable to the Kronecker definition 
except that the $b$--operation is replaced by  ``the  trivial 
star--operation'', called the    \it identity operation  \rm  
(for short, \it  $d$--operation\rm), acting as the identity map, i.e.
$$\Na(D)= \left\{ \frac{f}{g} \mid  f,g \in D[X],\ \boldsymbol{c}(f)^d
\subseteq \boldsymbol{c}(g)^d, \ \boldsymbol{c}(g) \mbox{ \rm  is 
invertible} \right\}.
$$  

A way to combine the ideas of the previous two paragraphs is to
view the e.a.b. property not as a property of a  
$\star$--operation,   as
has been done classically, but to view it as a property of a
certain class of ideals.  The e.a.b.  $\star$--operations   
should be
those for which every  nonzero   finitely generated ideal is an
``e.a.b.--ideal''  (Definition \ref{def:3.4}).   On the other 
hand, the invertible ideals should be
the only e.a.b.--ideals associated with the identity operation. So
given any  $\star$--operation,    we can combine the two 
definitions by
specifying that the content ideals of denominators must be
e.a.b.--ideals associated to the given  $\star$--operation.   

Note that,  from the beginning of the present paper,  we move from 
just  ``star--operations''  to the more general
setting of ``semistar--operations''  introduced by Okabe-Matsuda in 
1994 \cite{OM2}  (the definition is recalled in Section 2).   In 
fact, our generalizations work directly, and more naturally,  in this 
setting.

We also want to define our generalized function rings using the
mehod of intersecting Nagata rings of  quasilocal    overrings 
as is done
above.  Suppose then that we are
given a domain $D$ and a semistar operation  on $D$.  For
the $b$--operation we chose the class of all valuation overrings of
$D$ to define the Kronecker function ring and we chose the class of
all localizations of $D$ to define the Nagata ring of $D$.  When we
consider our discussion of e.a.b.--ideals above, we note that all
finitely generated ideals of $D$ extend to principal ideals in any
valuation overring.  On the other hand, invertible ideals are the 
only 
finitely generated ideals which extend
to principal ideals in every localization of $D$.  So the way to
proceed seems to be to combine the overring characterizations of the
Kronecker and Nagata rings by choosing the overrings in which the
e.a.b.--ideals extend to principal ideals.

\smallskip

Let $D$ be a domain and $\star$ a semistar operation
on $D$.  Given either a collection of overrings of $D$ or a ring 
of rational functions,    overring
of $D[X]$,  it is easy to define a new semistar operation on $D$ by
either extending to all of the overrings  of the collection   
and intersecting, or by
extending to the ring of rational functions and contracting back to
the quotient field of $D$.
We explore both of these mechanisms for
defining new semistar operations associated with  the given   
$\star$ and compare
the properties we obtain with those  proven   in the classical
Kronecker and Nagata settings.
  In particular, we deepen the study  of the ring of rational 
functions called  \it  the $\star$--Nagata ring  \rm  $\Na(D, \star)$ 
(Definition \ref{def:3.1}) and  \it  the $\star$--Kronecker function 
ring \rm  $\Kr(D, \star)$ (Definition \ref{def:3.2})
and we give a complete positive answer to the following question:\\
Is it possible to find a ``new''  integral domain of rational 
functions   
 denoted by $\KN(D, \star)$ \ (\it  `$`\star$--Kronecker-Nagata 
ring'',  \rm obtained  as  an 
intersection of  Nagata domains of quasilocal domains associated to   
a given arbitrary semistar operation $\star$ ) \ such that: 
\newline    $\bullet$ \;  $ 
\Na(D, 
\star) \subseteq \KN (D, \star) \subseteq  \Kr(D, \star)$\,;  
\newline  
$\bullet$ \;  $\KN(D, \star)$  ``generalizes'' at the same time 
$\Na(D, 
\star)$ and $\Kr(D, \star)$ and coincides with 
$\Na(D, \star) $ or $\Kr(D, \star)$, when the semistar 
operation  $\star$ assumes 
the ``extreme values'' of an interval $\star'  \leq \star 
\leq  \star''$\  (i.e. $\KN(D,\star') = \Na(D, \star')= \Na(D, \star) 
$ and $\KN(D,\star'') = \Kr(D, \star'')= \Kr(D, \star)$?

%%%%%%%%%%%%%%%%%%%%%%%%%%%%%%  SECTION 2
\section{Background}

In this section we give some definitions and some basic results, some
new and some not.

We begin by designating the following terms.

\begin{itemize}
\item $\boldsymbol{f}(D)$ is the set of all nonzero finitely 
generated fractional ideals
of $D$.
\item $\boldsymbol{F}(D)$ is the set of all nonzero fractional ideals 
of $D$.
\item $\overline{\boldsymbol{F}}(D)$ is the set of all nonzero $D$ 
submodules of $K$.
\end{itemize}

In 1994 A. Okabe and R. Matsuda  \cite{OM2}   introduced the 
notion
of a semistar operation.  A semistar operation is a map  $\,\star 
: 
\boldsymbol{\overline{F}}(D) \rightarrow
\boldsymbol{\overline{F}}(D)\,$,  $\,E \mapsto E^\star\,$   which 
obeys
the following axioms,   for
all $\,z \in K\,$, $\,z \not = 0\,$ and for all $\,E,F \in
\boldsymbol{\overline{F}}(D)$.

\vspace{.1in}

$\mathbf{ \bf (\star_1)} \; \; (zE)^\star = zE^\star \,; $

\vspace{.1in}

$\mathbf{ \bf (\star_2)} \; \;   E \subseteq F \;\Rightarrow\; E^\star
\subseteq
F^\star \,;$

\vspace{.1in}

$\mathbf{ \bf (\star_3)} \; \;   E \subseteq E^\star \; \; \textrm { 
and
}
\; \; E^{\star \star} := (E^\star)^\star = E^\star\,. $

\vspace{.1in}

The classical notion of a star operation  \cite{G}   involves a 
map from $\boldsymbol{F}(D)$
to $\boldsymbol{F}(D)$ which requires, in addition to the semistar 
axioms, that
 $(\alpha D)^{\star} = \alpha D$,    for each  nonzero    
principal ideal  $\alpha D$   of $D$.

The key difference here is that if $\star$ is a star operation then
$D^{\star} = D$, whereas $D^{\star}$ may be properly larger than $D$ 
 (possibly, $D^\star \in \boldsymbol{\overline{F}}(D) \setminus 
{\boldsymbol{F}}(D)$)   
if $\star$ is a semistar operation.  Note that if $\star$ is a
semistar operation on a domain $D$ then we obtain a classical star
operation on $D^{\star}$ when we restrict $\star$ to
$\boldsymbol{F}(D^{\star})$.

Now we give some basic information concerning definitions/terminology,
ge\-ne\-ral properties of semistar operations, and concerning the
construction of specific semistar operations on integral domains.

\begin{itemize}

\item As in the classical star operation setting, we associate to a 
semistar
ope\-ra\-tion $\,\star\,$ of $\,D\,$ a new semistar operation
$\,\star_f\,$ as follows.
 \rm Let $\,\star\,$ be a semistar operation of a domain $\,D\,$.  \ 
If
 $\,E \in \boldsymbol{\overline{F}}(D)\,$ we set:
$$
E^{\star_f} := \cup \{F^\star \;|\;\, F \subseteq E,\, F \in
\boldsymbol{f}(D)
\}\,.
$$
\noindent We call $\,\star_f\,$ \it the semistar operation of finite
type of
$D$ \rm associated to $\,\star\,.$ \ If $\,\star = \star_f\,$,\ we say
that
$\,\star\,$ is \it a semistar operation of finite type of $\,D\,$.  
\rm
\ Note that $\,\star_f \leq \star\,$ and $\,(\star_f)_f = \star_f\,$,\
so
$\,\star_f\,$ is a semistar operation of finite type of $\,D\,.$

\item \it A (semi)star operation $\star$ on $D$ \rm  is a semistar 
operation on
$D$ such that $D^{\star} = D$; i.e.  a semistar operation such
that: $$ \star\!\! \mid_{{\boldsymbol{F}}(D)}:
         {\boldsymbol{F}}(D) \rightarrow {\boldsymbol{F}}(D)$$
is a ``classical'' star operation  \cite[Section 32]{G}.  

\item $d_{D}$ denotes the identity (semi)star operation on $D$.

\item If $\iota : D \hookrightarrow T$ is the canonical embedding
of $D$ in the overring $T$ of $D$ and if $\star$ is a semistar 
operation
on $D$, then $\star_{\iota}$ is the semistar operation on $T$ defined,
for each $E \in \overline{\boldsymbol{F}}(T) \ (\subseteq
         \overline{\boldsymbol{F}}(D))$,
by
$$E^{\star_{\iota}} := E^\star\,.$$

\item If $T$ is an overring of $D$, we denote by $\star_{\{T\}}$
the semistar operation on $D$ defined as follows: for each
$E \in \overline{\boldsymbol{F}}(D)$:

         $$ E^{\star_{\{T\}}} := ET \,.$$
Obviously, $({\star_{\{T\}}})_{\iota} = d_{T}$\ (where $d_{T}$ denotes
the identity semistar operation on $T$)\,.

\item If $\{\star_{\lambda} \mid \lambda \in \Lambda \}$ is a
family of semistar operations on $D$ then
$\wedge \{\star_{\lambda} \mid \lambda \in \Lambda \}$ is the semistar
operation on $D$ defined as follows: for each
$E \in \overline{\boldsymbol{F}}(D)$:

   $$ E^{\wedge \{\star_{\lambda} \mid
         \lambda \in \Lambda \}} := \bigcap \{E^{\star_{\lambda}} \mid
         \lambda \in \Lambda \} \,.$$
In particular, if $\calbT:=\{T_\lambda \mid T_\lambda\in \Lambda\}$
is a given family of overrings of $D$, then $\wedge_{\calbT}$ denotes
the semistar operation
${\wedge \{\star_{\{T_\lambda\}} \mid \lambda \in \Lambda \}} $.

\item $b_{D}$ is the $b$--semistar operation on $D$, i.e.

   $$b_{D} := \wedge \{\star_{\{V\}} \mid
        \ V \mbox{ is a valuation overring of } D \}\,.$$

\end{itemize}

\bigskip

It seems natural in the context of the ``$\wedge$--construction''  
above,   i.e. $\wedge_{\calbT}$,  
to view semistar operations  ``extensions to the overrings'' , 
i.e. $\star_{\{T_\lambda\}}$,    as canonical components of
the semistar operation on $D$.  We have defined the
$b$--operation as a  $\wedge$--construction   using the 
valuation overrings.
In that setting we can think of the $b$--operation as being decomposed
into component semistar operations, each defined by extension to a
valuation overring.

A question that seems not to have been dealt with in the
literature (in either the star or semistar setting) is the extent to
which a given star (or semistar) operation on $D$ can be
approximated by one built from component parts  of the type   
``extensions to the overrings''    
of $D$ in the above manner.  One use for our generalized
Kronecker-Nagata theory is to associate  in a natural way    a  
semistar operation   defined by a $\wedge$--construction  
to a given semistar operation.

\bigskip

In the setting of star operations, the class of $\star$--ideals (i.e.
those ideals $I$ such that $I^{\star} = I$) assume a role of great
importance.  When $\star$ is a semistar operation there frequently are
no integral ideals of $D$ which are $\star$--ideals.  Instead we use
the following more general concept.

%% DEFINITION 2.1
\begin{deef} \label{def:2.1} \rm  Let $\,I \subseteq D\,$ be a 
nonzero ideal of $\,D\,$ and let
$\,\star\,$ be a semistar operation on $\,D\,$.  We say that $\,I\,$ 
is
\it a
quasi--$\star$--ideal of $\,D\,$ \rm if $\,I^\star \cap D = I\,$. \
Similarly, we
designate by \it  quasi--$\star$--prime \rm [respectively,\, \it
$\,\star$--prime
\rm] of $D$ a quasi--$\star$--ideal [respectively,\, an integral
$\,\star$--ideal]  of $\,D\,$  which is also a prime ideal. We 
designate by \it
 quasi--$\star$--maximal \rm
[respectively,\, \it  $\,\star$--maximal \rm]  of $\,D\,$  a maximal
element in the
set of all proper quasi--$\star$--ideals [respectively,\, integral
$\star$--ideals] of $\,D\,.$
\end{deef}

\vspace{.1in}

Note that if $\,I \subseteq D\,$ is a $\,\star$--ideal, it is also a
quasi--$\star$--ideal and, when $\,D = D^\star\,$
the notions of quasi--$\star$--ideal and
integral $\,\star$--ideal coincide.

We then give the following designations related to quasi-star ideals.

\begin{itemize}

\item  \; \quad $\QSpec^{\star}(D)$ is the set of all the 
quasi--$\star$--prime
ideals of $D$.

\item \; \quad $\calM({\star})$ is the set of all the maximal 
quasi--$\star$--ideals
of $D$.

\end{itemize}

 It is well known that \sl if $\star$ is a semistar operation of 
finite type then $\calM({\star})$  is nonempty  \cite[Lemma 2.3 
(1)]{FL3}.  \rm

\rm A particular important semistar operation (of finite type) on $D$ 
is the following:

\begin{itemize}

\item \; \quad
$
\widetilde{\star} := \wedge  \{\star_{\{D_Q\}} \mid
        \ Q \in \calM(\star_f) \} \ (= \wedge_{\calbS}\,,
       \mbox { where $\calbS :=\{ D_Q \mid
       \ Q \in \calM(\star_f) \}$)}\,.
       $
 \end{itemize}

       For the motivations,    examples    and the basic properties 
of this   type of    semistar 
operation cf. for instance  \cite[Corollary 2.7 and Remark 2.8]{FL3}.

%%%%%%%%%%%%%%%%%%%%%%%  SECTION 3
\section{{\rm e.a.b.}--ideals}

In this section we give the promised modification of the e.a.b.
condition motivated by the definitions of the Nagata rings and
Kronecker function rings.

We begin by giving more general definitions of the Kronecker function
ring and the Nagata ring.

 Recall that, for general semistar operations, we can consider the
notion of semistar invertible ideals.

%DEFINITION 3.1

\begin{deef}  \label{def:3.5} \rm  A fractional ideal
$I \in \boldsymbol{F}(D)$
is called \it  ${\star}$--invertible \rm if
$(II^{-1})^{\star} = D^{\star}$.
\end{deef}

For the motivations,   examples     and the basic properties of this 
type of 
invertibility see \cite{FP}.

%DEFINITION 3.2
\begin{deef}  \label{def:3.1}  \rm  Let $D$ be a domain   with 
quotient field $K$ and let $\star$ be
a semistar operation on $D$.  Set 
\[ \Na(D,\star) := \left\{\frac{f}{g}\mid f,g \in D[X], \co(f)^{\star}
\subseteq \co(g)^{\star},\ \co(g) \mbox{\rm\   is a 
$\star_{f}$--invertible ideal of } D \right\}.\]
This is called  \it  the $\star$--Nagata ring of $D$. \rm     If 
the semistar operation is the
identity operation,  i.e. $\star = d_D$,   then   
$\Na(D,d_D) $   coincides with the ``classical'' Nagata ring 
$D(X)$ of $D$.

\end{deef}

It is known that,  \sl for each $E  \in\overline{ \boldsymbol{F}}(D)$,
$$
E^{\widetilde{\star}}= E\Na(D, \star) \cap K\,,\;\; 
\mbox{\cite[Proposition 3.4 (3)]{FL3}}\,.
$$
  \rm

The definition of an e.a.b. semistar operation is as follows.

%DEFINITION 3.3
\begin{deef} \label{def:3.3} \rm Let $D$ be a domain and let $\star$
be a semistar operation on $D$.  Then we say that \it $\star$ is an
e.a.b. semistar operation \rm  provided $(IJ)^{\star}  \subseteq 
(IH)^{\star}$ 
implies $J^{\star} \subseteq H^{\star}$ whenever $I,J,H \in 
\boldsymbol{f}(D)$.  
(We say that
$\star$ is \it an a.b. semistar operation \rm if we weaken the 
hypotheses
to only require that $J,H$ lie in $\boldsymbol{\overline{F}}(D)$.)
\end{deef}

\smallskip

%REMARK 3.4
\begin{reem} \label {rk:3.4}    \bf (1) \rm  The paper \cite[in 
preparation]{FL6}   is devoted to  a deeper study on the relations 
between the e.a.b. and the a.b. semistar operations.    

  \bf (2) \rm  Recall that \it the e.a.b. semistar operation 
on $D$  of
finite type ${\star_{_{\! a}}}$ associated to a semistar
operation $\star$ \rm  is defined on $D$ by setting for
each $G \in \boldsymbol{f}(D)$:
$$
G^{\star_{_{\! a}}}:= \bigcup \left\{\ \left((GH)^\star : H^\star
\right) \mid\, H \in \boldsymbol{f}(D)\right\}\,
$$ 
\cite{J}, \cite{HK1} and \cite{HK2}.   
It is known that a semistar operation of finite type is
e.a.b. if and only if $\star = {\star_{_{\! a}}}$
(cf., for instance, \cite[Proposition 4.5]{FL1}).
\end{reem}

\bigskip 

\rm  Let $D$ be an integrally closed domain    with quotient field 
$K$ and let $\star$ be
 an e.a.b. semistar
operation on $D$.    Set
\[
 \Kr(D,\star) := \left\{\frac{f}{g} \mid  f,g \in D[X],\  g \neq 0, \ 
\co(f)\subseteq
\co(g)^{\star} \right\}\,. 
\]
\rm     
In the case where $\star$ is an e.a.b. \sl star \rm operation on  
(an integrally closed domain)   $D$ this
definition yields the  classical Krull's extension of the 
Kronecker function ring of $D$
associated with $\star$ \cite[Section 32 ]{G}.

\bigskip

This is not the most general definition of
the Kronecker function ring, but it is the one most suited to our 
program.   Recall that the reason customarily given for the 
assumption that
$\star$ be e.a.b. in the Krull's definition of the Kronecker function 
ring
is that the classical proof that $\Kr(D)$ is a ring does not
work unless $\star$ is assumed to be e.a.b..  The use that the
e.a.b. property is put to in this context is to show that an
equation of the form $(c(g)J)^{\star} = (c(g)H)^{\star}$ implies
that $J^{\star} = H^{\star}$ where $J,H \in \boldsymbol{f}(D)$ and
$g$ is a nonzero polynomial in $D[X]$ and the denominator of a 
rational
function $\frac{f}{g}$ in $\Kr(D,\star)$.

\medskip

A general definition for the Kronecker function ring, without 
restrictions on $D$ and $\star$, is recalled next.

%DEFINITION 3.5
\begin{deef} \label{def:3.2} \rm  Let $D$ be any integral domain (not 
necessarily integrally closed)   with quotient field $K$ and let 
$\star$ be
 any semistar
operation on $D$ (not necessarily e.a.b.).    Set
\[
\Kr(D,\star) := \left\{ \frac{f}{g} \ \mid 
 \begin{array}{rl}   & \hskip -10 pt f,g \in D[X],\ g\neq 0,\ \mbox{ 
such that there exists  } h  \in D[X], \\
  &  h\neq 0,\ \mbox{ with   }  \co(fh) \subseteq
\co(gh)^{\star}
\end{array}
\right\}\,.
\]
This is called  \it  the $\star$--Kronecker function ring of $D$, 
\cite[Theorem 5.1]{FL1}. \rm 
Obviously if  $\star$ is an e.a.b. semistar
operation on  $D$, then this definition coincides with the previous 
one.
\end{deef}

\smallskip

 It is known that,  \sl for each $E  \in \overline{ 
\boldsymbol{F}}(D)$,
$$
E^{\star_a}= E\Kr(D, \star) \cap K\,,\;\; \mbox{\cite[Proposition 4.1 
(5)]{FL3}}\,.
$$
   \rm

Since invertible ideals can be cancelled in any conceivable context,
it is clear that the ``modified'' e.a.b.   property,   that we 
want to introduce for generalizing semistar Nagata rings and 
Kronecker function rings,   should be a   ``cancellation--type'' 
  property.

The definition is actually quite straightforward.

%DEFINITION 3.6
\begin{deef} \label{def:3.4} \rm Let $D$ be an integral domain  and 
let
$\star$ be a semistar operation on $D$. If $F$ is in
$\boldsymbol{f}(D)$, we say that \it $F$ is a $\star$--e.a.b--ideal 
\rm  if
$(FG)^{\star} \subseteq (FH)^{\star}$, with $G,\ H \in
\boldsymbol{f}(D)$, implies that $G^{\star}\subseteq H^{\star}$.
As with the semistar operations, we say that \it $F$ is an 
a.b.--ideal \rm 
if the conclusion holds with the requirement weakened to say that
$G,\ H \in  \boldsymbol{\overline{F}}(D)$.  
\end{deef}

It is clear that \sl a semistar operation $\star$ on a domain $D$ is
e.a.b. if and only if every finitely generated ideal of $D$ is a
$\star$--e.a.b. ideal. \rm

%REMARK 3.7
\begin{reem} \label{rk:3.5}
\rm It is clear that  \sl invertible ideals are $\star$-e.a.b. for any
semistar operation $\star$. \rm   In fact, if the semistar operation 
in
question is the identity operation $d$, the $d$--e.a.b. ideals of
a domain $D$ correspond to what D.D.  Anderson and D.F. Anderson 
called
quasi-cancellation ideals.  They proved  that,  in the finitely 
generated setting,    the
$d$--e.a.b. ideals are exactly the invertible ideals \cite[Lemma 1
and Theorem 1]{AA}.
\end{reem}

It is easy to see that,  in general, \sl a finitely generated  
$\star$--invertible ideal is also
$\star$--e.a.b..    \rm  Unlike the case for the identity operation
though  (Remark \ref{rk:3.5}),   it is not true in general that 
all   (possibly, finitely generated)   $\star$--e.a.b. ideals
are $\star$--invertible.  
For instance, let
$D$ be a Noetherian domain of dimension greater than one and let
$M$ be a maximal ideal  of $D$   of height greater than one.  
 Since the $b$--operation  is an e.a.b. star operation on $D$, 
then  $M$ (which is finitely generated) is a
$b$--e.a.b. ideal, by the observation preceding Remark \ref{rk:3.5}. 
  However, since $M$ has height greater than one,  it is not   an
invertible ideal of $D$.   Hence $MM^{-1} = M$.  Then 
$(MM^{-1})^{b} = M^{b} =
M$.  Hence, $M$ is not $b$--invertible.

\bigskip 

We close this section with a collection of basic results concerning
$\star$--e.a.b. ideals, invertible ideals, and $\star$--invertible
ideals.

It is known that \sl if $\star$ is an e.a.b. star operation on an
integral domain $D$, then there exists an a.b. star operation
$\ast$ on $D$ such that $ \star\!\!\mid_{\boldsymbol{f}(D)} =
\ast\!\!\mid_{\boldsymbol{f}(D)}$ \rm  \cite[Corollary 32.13]{G}.  
  This
motivates our next statement proven in \cite{FL6}.

%LEMMA 3.8
\begin{leem} \label{le:1.5}   \sl Let $D$ be an integral domain 
and let $\star$ be a
semistar operation on $D$. 

\begin{enumerate}

\item[\bf (1)] \sl If $\star = \star_{_{\! f}}$, then: 
$$ \star \mbox{ \sl  is an e.a.b. semistar operation if and only if 
 $\star$ is an a.b. semistar operation}.$$

\it\item[\bf (2)] \sl Let  $F \in \boldsymbol{f}(D)$, then:
$$ \mbox{F is $\star$--e.a.b. if and only if  F is $\star_{_{\! 
f}}$--a.b.} $$
\end{enumerate}

\end{leem}  \vskip - 18   pt \hfill $\Box$

\medskip

The following result is known \cite[Theorem 2.23]{FP}.

% LEMMA 3.9
\begin{leem} \label{lm:1.3} \sl Let $\star$ be a semistar operation
on an integral
domain
$D$.  Let $F \in \boldsymbol{f}(D) $,
then the following are equivalent:

\brom  % 1
 \bf \item \sl $F$ is
${\star_{_{\!f}}}$--invertible;
 %2
\bf \item \sl $FD_Q$  is invertible as a fractional ideal of  
$D_Q$,   for each
$Q\in\calM({\star_{_{\!f}}})$;
 %3
\bf \item \sl $F\Na(D,\star)$   is invertible as a fractional 
ideal of $\Na(D,\star)$.     \erom
\vskip -14pt \hfill $\square$
\end{leem}

% REMARK 3.10
\begin{reem} \label{rk:1.2} \rm

  \bf (1) \rm  Let $F \in \boldsymbol{f}(D) $.  As a consequence of
Lemma \ref{lm:1.3}, note that, since $\calM(\star_{_{\!f}})=
\calM(\tilde{\star})$ \cite[Corollary 3.5]{FL3} (2), then:
             \sl

             \centerline{ $F$ is $\star_{_{\!f}}$--invertible if and
only if $F$ is
     $\tilde{\star}$--invertible,} \rm

              \noindent  cf. \cite[Proposition 2.18]{FP}.

 \bf (2) \rm Let $F \in \boldsymbol{f}(D)
$.    From \cite{FL6} recall that \ \sl $F$ is $\star$--e.a.b.  
[respectively, $\star$--a.b.]   if 
and only if  $\left(
               (FH)^\star : F^\star \right) = H^\star$, for each $H
\in
               \boldsymbol{f}(D)$ [respectively, for each $H \in
              \overline{ \boldsymbol{F}}(D)$]. \rm  \
(Note that  $\left((FH)^\star : F^\star \right) = \left(
               (FH)^\star : F\right) $, and so the previous
equivalences can be stated in a formally slightly different way.)   
   \end{reem}

%%%%%%%%%%%%%%%%%%%%%%%%%%%   SECTION  4

\section{Some distinguished classes of overrings}

We begin   by    considering a class of overrings of a domain 
$D$ associated
with a semistar operation $\star$ which have already been well 
studied.

%DEFINITION 4.1 
\begin{deef}  \label{def:4.1}  \rm   Let $\star$ be a semistar 
operation
on an integral
domain $D$.  We say that an overring $T$ of $D$ is a
{\it $\star$--overring } of $D$ provided for each
$F \in \boldsymbol{f}(D) $ we have $F^{\star} \subseteq FT$ (or
equivalently $F^{\star}T = FT$).
\end{deef}

The following lemma gives some basic results concerning
$\star$--overrings.

        %LEMMA 4.2
\begin{leem} \label{lm:1.1} \sl Let $D$ be an integral domain
with quotient field $K$ and let $\star$ be a semistar operation on
$D$.\bara

% 1
\bf \item [(1)]\sl  The following are equivalent:
\begin{enumerate}
\bf \item[(i)] \sl
$T\ \mbox{ is a $\star$--overring of }\ D $\,;

\bf \item[(ii)] \sl $ T\ \mbox{ is a ${\star_{_{\!f}}}$--overring
of }\ D $\,;

\bf \item[(iii)] \sl  $\star_{_{\!f}}\leq \star_{\{T\}}\,,$ \;
(i.e. $E^{\star_{_{\!f}}} \subseteq ET\,,\; \forall\
E\in \overline{\boldsymbol{F}}(D)\,). $

\bf \item[(iv)] \sl  $({\star_{_{\!f}}})_{\iota} =d_{T}$\,.
\ealf In particular, if $T$ is a $\star$--overring of  $D$, then
$D^\star \subseteq T^{\star_{_{\!f}}} = T$.

% 2
\bf \item[(2)] \sl  Every overring of a $\star$--overring is a
$\star$--overring.
% 3
\bf \item[(3)] \sl $K$ is a $\star$--overring of  $D$, for
each semistar operation $\star$ on $D$.
%4
\bf \item[(4)] \sl $D^\star$ is a $\star$--overring of $D$ if
and only if $\star_{_{\!f}} = \star_{\{D^\star\}}$. More
generally, $T$ is  a $\star$--overring of  $D$ and $T =
D^\star$ if and only if $\star_{_{\!f}} = \star_{\{T\}}$.
%5
\bf \item[(5)] \sl  If $\star_{1} \leq \star_{2}$ are two
semistar operations on $D$, then:
$$ T \ \mbox{ is a $\star_{2}$--overring of } \ D \;\;
\Rightarrow \;\;  T \  \mbox{ is a $\star_{1}$--overring
of }\ D\,. $$
%6
             \bf \item[(6)] \sl A B\'ezout overring $B$ of $D$ is a
$\star$--overring  of $D$ if and only if $ B = B^{\star_{_{\!f}}}$ \.
In
             particular,  a valuation overring $V$ of $D$ is a
$\star$--overring  of $D$ if and only if $ V = V^{\star_{_{\!f}}}$\
(in this
situation,  $V$ is
             called a \rm $\star$--valuation overring of $D$\sl ).
%7
              \bf \item[(7)] \sl The valuation overrings of a
$\star$--overring $T$ of an
              integral domain $D$ coincide with the
$\star$--valuation
overrings of
              $D$ containing $T$.
% 8
              \bf \item[(8)] \sl Let $T$ be a $\star$--overring of
$D$ and let $\iota: D
              \hookrightarrow T$ be the canonical inclusion.  Then:
              $$
              \Kr(D, \star_{\{T\}})=
             \Kr(T, d_{T}) = \Kr(T, (\star_{_{\!f}})_{\iota}) =
\Kr(T, b_{T})\,.
              $$
%9
              \bf \item[(9)] \sl If $N$ is a prime ideal of a
$\star$--overring $T$ of $D$
              and if $N \cap D \neq 0$, then $N \cap D$ is a
              quasi--$\star_{_{\!f}}$--prime of $D$.
%10
             \bf \item[(10)] \sl Let $\calT := \{T_\lambda \mid
\lambda \in \Lambda\}$ be a family of overrings of $D$.
          The semistar operation $\wedge_{\calT}$ (defined  in 
Section 2)    is
such that:  $$
          E{T_\lambda}  = E^{ \wedge_{\calT}}T_\lambda =
(ET_\lambda)^{ \wedge_{\calT}}\,,
          $$
  for each  $E\in \overline{\boldsymbol{F}}(D)$.  In
particular, each $T_\lambda $ is a $\wedge_{\calT}$--overring of
$D$.  \eara \end{leem}

             \bf{Proof.} \rm \bf (1) \rm (i) $\Leftrightarrow$ (ii)
$\Leftrightarrow$
             (iii) are obvious consequences of the definition. \\ 
(iii)
$\Rightarrow$
             (iv): $d_{T }\leq (\star_{_{\!f}})_{\iota}\leq
             (\star_{\{T\}})_{\iota} =d_{T}$.  \\   (iv) $\Rightarrow$
             (iii): For each $E \in \overline{\boldsymbol{F}}(D)$, $
      ET= (ET)^{(\star_{_{\!f}})_{\iota}} =(ET)^{\star_{_{\!f}}}
      \supseteq  E^{\star_{_{\!f}}}$.

 \bf (2) \rm is an
easy consequence of (1).

\bf (3) \rm  and \bf (5) \rm are obvious.

\bf (4) \rm   follows from (1) and from the fact that, in general,
$\star_{\{D^\star\}} \leq \star_{_{\!f}}$.

             \bf (6) \rm The ``only if'' part is obvious from (1). 
             For the ``if'' part recall that,
for each $F
             \in \boldsymbol{f}(D)$, there exists a nonzero element
$x
\in K$ such
             that $FB = xB$, thus $ F^\star \subseteq
(FB)^{\star_{_{\!f}}} =
             (xB)^{\star_{_{\!f}}}=
             xB^{\star_{_{\!f}}}= xB = FB$.

             \bf (7) \rm follows from (2) and (6).

             \bf (8) \rm Since $T$ is a $\star$--overring of $D$,
             then,  by (1), \  $ (\star_{_{\!f}})_{\iota}= d_{T}$.
              Therefore, by (7),
             $V$ is a
($d_{T}$--)valuation overring of $T$ if and only if $V$ is a
$\star$--valuation overring of $D$ containing $T$ if and only if
$V$ is a
$\star_{\{T\}}$--valuation overring of $D$.
The
conclusion is a straightforward consequence of the fact that, if
$\ast$ is a semistar operation on an integral domain $R$, then
$\Kr(R,
\ast) = \cap \{W(X) \mid W  \mbox{ is a $\ast$--valuation overring of
}
R\}$ \cite[Theorem 3.5]{FL2}.

              \bf (9) \rm $(N \cap D)^{\star_{_{\!f}}} \subseteq (N
\cap
D)T \subseteq
              N$, hence $N \cap D \subseteq (N \cap
D)^{\star_{_{\!f}}}
\cap D \subseteq
              N\cap D$.

              \bf (10) \rm  Note that
$(ET_\lambda)^{ \wedge_{\calT}} = \bigcap \{ET_\lambda T_\mu \mid \mu
\in \Lambda\} = ET_\lambda \cap (\bigcap  \{ET_\lambda T_\mu \mid \mu
\in \Lambda\,, \; \mu \neq \lambda\})  = ET_\lambda $.
             \hfill $\Box$

%COROLLARY 4.3
 \begin{coor} \label{co:1.5}  \sl Let $\star$ be a semistar operation 
on an integral domain
$D$, let $F \in \boldsymbol{f}(D) $ be
${\star_{_{\!f}}}$--invertible and let $(L, N)$ be a local
$\star$--overring of $D$.  Then $F\! L$ is a principal fractional
ideal of $L$.
\end{coor}

\bf Proof.  \rm Recall that an invertible ideal in a local domain
is principal \cite[Theorem 59]{Ka}.  By   {Lemma \ref{lm:1.3}},  
we know that $F\! D_{Q}$ is principal in $D_{Q}$, for each $Q \in
\calM({\star_{_{\!f}}})$.  Hence, $F\! D_{N\cap D}$ is principal
in $D_{N\cap D}$, because $N\cap D$ is a
quasi--${\star_{_{\!f}}}$--prime ideal of $D$ (Lemma \ref{lm:1.1}
(9)), hence $N\cap D \subseteq Q$, for some $Q \in
\calM({\star_{_{\!f}}})$.  The conclusion follows immediately,
since $(L, N)$ dominates $(D_{N\cap D}, {N\cap D})$.   \hfill
$\square$

 % REMARK 4.4
             \begin{reem}  \label{rk:4.3}  \rm Let $D$ be an integral 
domain
        and let
         $\star$ be a semistar operation on $D$.
         
         \bf (1)   \rm  An overring $T$ of
$D$
        such that $T = T^{\star_{_{\!f}}}$ is not necessarily a
        $\star$--overring of $D$.

        For instance, let $\calF$ be a
        localizing system of $D$, and let $\star := \star_{\calF}$ be
        the semistar operation on $D$,  defined by $
        E^{\star_{\calF}} := E_{{\calF}} := \cup \{(E:I) \mid I \in
        \calF\}$, for each $E \in \overline{\boldsymbol{F}}(D)$ (cf.
        \cite[Proposition 2.4]{FH}).  Set
        $T := D^{\star_{\calF}} =D_{{\calF}}$, then
        in general, $ED_{{\calF}} \subsetneq E_{{\calF}}$.  More
        precisely,  $\star_{\{D^{\star_{\calF}} \}} ={\star_{\calF}}$
        if and only if $D^{\star_{\calF}}$ is $D$--flat  and $\calF
        =\{ I \mbox{ ideal of } D \mid ID^{\star_{\calF}} =
        D^{\star_{\calF}}\}$ (cf. \cite[Proposition 2.6]{FH}).

        \bf (2) \rm We have mentioned in the proof of Lemma 
\ref{lm:1.1} (8) that:
      $$
      \Kr(D, \star) = \bigcap\{ V(X) \mid V \mbox{  is a 
$\star$--valuation overring of  }  D  \} \quad  \mbox{\cite[Theorem 
3.5]{FL2}.}
       $$
      From this property and from the observation following 
Definition \ref{def:3.2} it is possible to prove that:
     $$
     \star_a = \wedge_{\calV(D, \star)}\,, \; \mbox{ where  } 
\calV(D, \star):= \{ V\mid V \mbox{ is a $\star$--valuation overring 
of  } D \}\,.
     $$ 
 \end{reem}

In the introduction we alluded to classes of quasilocal overrings of
a domain $D$ that are associated with a given semistar operation on 
$D$.
The two classes of quasilocal overrings that arise in the Kronecker 
and
Nagata settings are localizations and valuation overrings.  Note
first that a finitely generated ideal of a domain $D$ is invertible
if and only if it is locally principal.  Also note that every finitely
generated ideal of $D$ extends to a principal ideal in any valuation
overring of $D$.  Since the collection of ideals we have been
concerned with are the $\star$--e.a.b. ideals it seems reasonable 
that what
we need are quasilocal overrings of $D$ in which each $\star$--e.a.b.
ideal extends to a principal ideal.  It turns out that assuming just
this property is not quite sufficient.  We give two refinements, each
of which we will use to create separate generalized function rings.

%DEFINITION 4.5
\begin{deef}  \label{def:4.4}  \rm Let $\star$ be a semistar 
operation on an integral
domain $D$ and let $T$ be a quasilocal overring of $D$. 

We say that
$T$ is \it  a $\star$--monolocality of $D$  \rm provided 
$T^{\star_{f}} = T$ and
every $\star$--e.a.b. ideal of $D$ extends to a principal ideal in 
$T$. 

We say that
$T$ is \it a strong $\star$--monolocality of $D$ \rm  provided $T$ is 
a
$\star$--overring of $D$ and
every $\star$--e.a.b. ideal of $D$ extends to a principal ideal in 
$T$.
\end{deef}

\smallskip

Let $\star$ be a semistar operation on a domain $D$ and let
$T$ be a $\star$--overring of
$D$. Then $T^{\star_{f}} = T$.  It follows then that \sl  a
strong $\star$--monolocality is a $\star$--monolocality. \rm

Also note that if $\star$ is the identity semistar operation on a
domain $D$ then every quasilocal overring of $D$ is a strong
$\star$--monolocality (and hence they are all 
$\star$--monolocalities).
At the opposite extreme, if $\star$ is the $b$--operation on a domain 
$D$,
then the collection of all strong $\star$--monolocalities and the
collection of all $\star$--monolocalities are both equal to the
collection of all valuation overrings of $D$.

Our goal for this
concept then is to identify semistar operations such that the
collection of $\star$--monolocalities (strong or not) is not all
quasilocal overrings and does not consist entirely of valuation
domains.  Or, conversely, to identify collections of overrings and
associate semistar operations,   using the 
$\wedge$--constructions,   which will give the collection of
overrings back as (strong) $\star$--monolocalities.

\medskip

We adopt the following notation.

Set: $$ \begin{array}{rl}  \calL:= & \hskip -5pt
\calL(\star):=  \calL(D, \star):=\{ L\mid L \mbox{ is a } \star\mbox
{--monolocality of } D \}\,,\\ {\calL}' := & \hskip -5pt  {\calL}'
(\star) := {\calL}' (D, \star) :=  \{ L' \mid L' \mbox{ is a
strong--$\star$--monolocality of } D \}\,.\\
\end{array}
$$

 Note that the set:
$
\calV:= \calV(\star):= \calV(D, \star) \ (= \{ V \mid V  \mbox{ is a
$\star$--valuation overring of } D \})
$
is obviously a subset of ${\calL}'$.

% LEMMA 4.6
\begin{leem} \label{lm:1.4} \sl Let $\star$ be a semistar operation 
on an
integral domain $D$. \bara
              %1
\bf \item \sl  $\calL(\star) =  \calL(\star_{_{\!f}})$
 \ and \
$ {\calL'}(\star)  =  {\calL'}(\star_{_{\!f}})$.  %2
\bf \item \sl A  quasilocal   overring $T$ of a 
$\star$--monolocality
of $D$ is also a $\star$--monolocality  of $D$ if and only if
$T = T^{\star_{_{\!f}}}$.  %3
\bf \item \sl  A  quasilocal   overring $S$ of a
strong--$\star$--monolocality of $D$ is always a
strong--$\star$--monolocality  of $D$.
   %4
\bf \item \sl  Let $F := (a_{1},\ a_{2},\
\ldots, a_{n})\!  D \in \boldsymbol{f}(D) $ be ${\star}$--e.a.b.
and let $L$ be a $\star$--monolocality of $D$.  Then
$$
 F\! L= F^\star\!  L = (F\!  L)^{\star_{_{\!f}}} =
a_{i}L\,, $$

for some $i$, with $1\leq i \leq n$.
\eara
\end{leem}

\bf Proof.  \rm \bf (1) \rm is obvious and \bf (2) \rm    and  \bf
(3) \rm follow from Lemma
\ref{lm:1.1} (1) and (2).

\bf (4) \rm   We start by recalling the following well known fact:

\bf Claim. \sl Let $F =(a_{1},\ a_{2}, \ldots,\ a_{n})D \in
\boldsymbol{f}(D)$ and let $L$ be a  quasilocal   overring of 
$D$. If $F\! L$
is principal in $L$, then, for some $i$, with $1 \leq i \leq
n$, $F\!L= a_{i}L$. \rm

If $F\!  L= (a_{1},\ a_{2}, \ldots,\ a_{n})
L= zL$, for some $z \in L$ then,  for each $i$, with
$1 \leq i \leq n$,
we can find a nonzero $x_{i} \in L$, such that $x_{i} z = a_{i}$.
Therefore, $zL=(a_{1},\ a_{2}, \ldots,\ a_{n})\!  L= (x_{1},\ x_{2},
\ldots,\ x_{n}) L \cdot zL$, hence $(x_{1}, x_{2},\ldots, x_{n})
L= L$.  Being $L$  quasilocal   then, for some $i$, with $1 \leq 
i \leq n$, we
have that $x_{i}$ is a unit in $L$, thus $zL= x_{i}zL=a_{i}L$.

Now we conclude the proof  of (4). Since
 $F$ is ${\star}$--e.a.b. and
$L \in \calL$, then $F\!L$ is principal, and thus, for
some  $i$, with $1 \leq i \leq n$, $(F\!  L)^{\star_{_{\!f}}} =
(a_{i}L)^{\star_{_{\!f}}}=
a_{i}{L}^{\star_{_{\!f}}} = a_{i}L$ (being $L = {L}^{\star_{_{\!f}}}$
). Therefore $a_{i}L\subseteq  F\!L\subseteq
F^{\star}L\subseteq
(F\!  L)^{\star_{_{\!f}}} = (a_{i}L)^{\star_{_{\!f}}}=
 a_{i}L$.  \hfill $\Box$
 \vskip 12pt

%%%%%%%%%%%%% SECTIION 5
\section{Generalized Kronecker--Nagata rings}

Now we turn to the construction of the generalized Kronecker and
Nagata rings.  We define two classes of rings which we refer to
as \it  Kronecker-Nagata ring \rm   (for short,   KN)  and  
\it   Strong Kronecker-Nagata ring \rm   (for short,   $\SKN$)  
according to whether we
use monolocalities or strong monolocalities.

 %PROPOSITION  5.1
             \begin{prro} \label{lm:1.7} \sl      Let $\star$
 be a semistar operation on an
             integral domain $D$.   Set:
             $$ \KN(D, \star) := \cap \{ L(X) \mid L \in \calL(D,
\star)\}\,, \;\;\;\;\;  \SKN(D, \star) :=\cap \{ L'(X) \mid L' \in
\calL' (D, \star)\}\,, $$
then:
             $$  \begin{array}{rl}
             \Na(D, \star) \subseteq & \hskip -5pt  \SKN(D, \star)
\subseteq \Kr(D,
\star)\,,\\ \KN(D, \star) \subseteq & \hskip -5pt \SKN(D, \star)\,.
\end{array} $$   \end{prro}
               \bf Proof.  \rm  By  Lemma
             \ref{lm:1.1} (9) and \cite[Proposition 3.1 (4)]{FL3},
we know that $ \Na(D, \star) = \cap \{D_{Q}(X) \mid Q
\in
             \QSpec^{\star_{_{\!f}}}(D) \} \subseteq  $ $\cap
\{D_{N'\cap
D}(X) $ $\mid (L', N') \in
             \calL' \}$ $ \subseteq \SKN(D, \star)$.  The 
inclusions $\KN(D, \star)  \subseteq \SKN(D,
\star)
\subseteq $ $\Kr(D,
             \star)$
follow
from the fact that $\calV(D, \star) \subseteq \calL' (D, \star)  
\subseteq \calL(D, \star) $.     \hfill $\Box$

\vskip 12pt  We have shown that the  Strong 
$\star$--Kronecker-Nagata ring   
$\SKN(D,\star)$ lies properly in between the   $\star$--Nagata 
   ring and the
 $\star$--Kronecker   function ring.  We have also shown that the
 $\star$--Kronecker-Nagata ring    $\KN(D,\star)$ lies inside 
the  $\star$--Kronecker   
function ring.  We will show later (Theorem \ref{th:1.21} (7))
that, in general, $\Na(D, \star) \subseteq  \KN(D, \star)$.
\vskip 12pt

Proposition \ref{lm:1.7} gives a positive result concerning
containment relations
for $\KN$ and $\SKN$.  The containment/inequality relations between
these concepts is not always as clean as we would like however.  For
example, it seems reasonable that if $\star_{1} \leq \star_{2}$ are
semistar operations on a domain $D$ then we would have
$\KN(D,\star_{1}) \subseteq \KN(D,\star_{2})$ and
$\SKN(D,\star_{1}) \subseteq \SKN(D,\star_{2})$.  
In Example \ref{ex:1.26} we give an example of a star operation 
$\star$ on a two-dimensional   Noetherian   local   integrally 
closed   domain $D$ such that 
$b < \star$ and yet $\KN(D,\star) = \Na(D, \star)   =D(Z) \subsetneq  
\KN(D, b) =\Kr(D, b) $.     Hence, in general, 
we do not get the containment we wish for $\KN$.  We do not know 
whether $\SKN$ behaves well with regard to containment and inequality 
or not.  We can give a positive result in this
direction when $\star_{2}$ is a stable semistar operation.

Recall
that  \it a semistar operation $\star$ \rm is \it stable on $D$ \rm 
provided
$$
(E \cap F)^\star = E^\star \cap F^\star,  \;\, \textrm {  for all} \;
E,F \in \boldsymbol{\overline{F}}(D) \,.
$$

  %COROLLARY 5.2
               \begin{coor} \label{cor:1.8}  \sl Let $\star_{1}
 \leq \star_{2}$ be two semistar operations on an
             integral domain $D$.    For $i =1, 2$, set:   $$
\calL_i := \calL(D, \star_i)\,,\;\;\;\;\;\;   {\calL'}\!_i :=
{\calL'}(D, \star_i)\,.   $$  Assume that $\star_{2}$ is
             stable.
              \bara
             % 1
             \bf \item \sl  Let $F \in
             \boldsymbol{f}(D)$. If $F$ is $\star_{1}$--e.a.b. then
             $F$ is also $\star_{2}$--e.a.b..
               % 2
              \bf \item \sl  ${\calL}_1 \supseteq {\calL}_2$ \ and \
${\calL'}\!_1 \supseteq {\calL'}\!_2
             $.

         % 3
               \bf \item \sl  $\KN(D, \star_{1}) \subseteq \KN(D,
              \star_{2})$ \ and \ $\SKN(D, \star_{1}) \subseteq
\SKN(D,
              \star_{2})$.
              \eara
              \end{coor}

              \bf Proof.  \rm \bf (1) \rm  is a consequence of
               Remark \ref{rk:1.2} (2),   since:
              $$
      \left(
               (FH)^{\star_{1}} : F\right) = H^{\star_{1}}  \;
               \Rightarrow\left(
               (FH)^{\star_{1}} : F\right)^{\star_{2}}=
               \left(H^{\star_{1}}\right)^{\star_{2}}\,,$$
             therefore, by the stability of $\star_{2}$ \cite[Theorem
2.10 (B)]{FH}:
              $$ H^{\star_{2}} \subseteq  \left(
               (FH)^{\star_{2}} : F\right)  = \left(\left(
               (FH)^{\star_{1}}\right)^{\star_{2}} : F\right)
               =  \left(\left(
               (FH)^{\star_{1}} : F\right)\right)^{\star_{2}}=
               H^{\star_{2}}\,,
            $$
               for each $H \in
               \boldsymbol{f}(D)$.

               \bf (2) \rm  follows from (1), from Lemma
               \ref{lm:1.1} (5) and from the fact that, if $T$ is a
 quasilocal   overring of
              $D$ such that $T= T^{(\star_{2})_{_{\!f}}}$,  then
              necessarily $T= T^{(\star_{1})_{_{\!f}}}$.  
              
              \bf (3)
\rm is a
               trivial consequence of (2). \hfill $\Box$

\bigskip

We have noted that the $\KN$ and $\SKN$ constructions are not always
well behaved in terms of preserving relationships between distinct
semistar operations.  Nonetheless, it seems worthwhile to pursue this
idea with regards to the operations  $\star_{a}$   and   
$\widetilde{\star}$  
associated to the semistar operations $\star$ on a domain $D$.  The
semistar operations $\star_{a}$ and $\widetilde{\star}$ are generally 
well behaved and
the results work out as we would hope.

We need a preparatory lemma first.

   %PROPOSITION  5.3
              \begin{prro} \label{prop:1.12} \sl  Let $\star$ be a
semistar operation on an
             integral domain $D$. Then:
  \bara

              % 1
             \bf \item \sl
For each $Q \in
\calM({\star_{_{\!f}}})$, $D_{Q}$ is
strong--$\widetilde{\star}$--monolocality of $D$.

% 2
\bf \item
\sl If $F \in
             \boldsymbol{f}(D) $, then $F$ is
             ${\widetilde{\star}}$--e.a.b. if and only if $F$ is
             ${\widetilde{\star}}$--invertible. \eara
\end{prro}
\bf Proof. \rm \bf (1) \rm  It is clear from the
definition of the
              semistar operation $\widetilde{\star}$ that, for each
$Q \in
\calM({\star_{_{\!f}}}) \ (= \calM(\widetilde{\star})$,
\cite[Corollary 3.5]{FL3}), $D_{Q}$ is a  quasilocal  
$\widetilde{\star}$--overring of $D$ since,  for each $F \in
             \boldsymbol{f}(D) $, $F^{\widetilde{\star}}\subseteq
             FD_{Q}$.

             Note, more generally, that for each $Q \in
\calM({\star_{_{\!f}}})$ and for each $E \in
             \overline{\boldsymbol{F}}(D) $:
             $$ED_{Q} = E^{\widetilde{\star}}D_{Q}=
             (ED_{Q})^{\widetilde{\star}}\,,$$
 (Lemma \ref{lm:1.1} (10), since $\widetilde{\star} =
\wedge_{\boldsymbol{\calS}}$, with ${\boldsymbol{\calS}} := \{ D_Q
\mid Q \in  \calM({\star_{_{\!f}}}\}$).

\rm \bf (2) \rm Let $F\in \boldsymbol{f}(D)$ be a
${\widetilde{\star}}$--e.a.b.,  thus
$\left((FH)^{\widetilde{\star}} : F\right) =
   H^{\widetilde{\star}}$ and so:
   $$
    HD_{Q} =H^{\widetilde{\star}}D_{Q} =
\left((FH)^{\widetilde{\star}} : F\right)D_{Q}
    =
   \left((FH)^{\widetilde{\star}}D_{Q} : FD_{Q}\right) =
   \left(FHD_{Q} : FD_{Q}\right)\,,
  $$
for each $H \in
\boldsymbol{f}(D)$,  i.e. $FD_{Q}$ is a quasi--cancellation ideal of
$D_{Q}$  or,
equivalently, it is a principal fractional ideal of $D_{Q}$, for each
$Q \in
\calM({\star_{_{\!f}}})$ (Remark \ref{rk:3.5}).

             From Lemma
             \ref{lm:1.3}  and Remarks \ref{rk:1.2} (1) and  
\ref{rk:3.5},    we deduce that   if $F \in
             \boldsymbol{f}(D) $, then $F$ is
             ${\widetilde{\star}}$--e.a.b. if and only if $F$ is
             ${\widetilde{\star}}$--invertible.           \hfill
$\Box$  \vskip 12pt

%PROPOSITION 5.4
\begin{prro} \label{pr:1.12} \sl
Let $\star$ be a semistar operation on an
             integral domain $D$, then:
               \bara
                 % 1

\bf \item \sl $\Na(D, \star) = \Na(D, \widetilde{\star}) = \SKN(D,
             \widetilde{\star})\,.$
             %2
             \bf \item \sl  $
\SKN(D,  \star_{a}) =\Kr(D, \star_{a}) = \Kr(D, \star)\,.$
\eara     \end{prro}
\bf Proof. \rm  \bf (1) \rm  By
Proposition \ref{prop:1.12} we
             know that   $ \calL'(D, \widetilde{\star})
\supseteq      \{ D_{Q} \mid Q \in
\calM({\star_{_{\!f}}})\}$ and if $(L', N')  \in \calL'(D,
\widetilde{\star})$ then $ L' \supseteq D_{N'\cap D} \supseteq D_Q$,
where $Q$ is any prime ideal in $\calM({\star_{_{\!f}}}) \
(=\calM({\widetilde{\star}}))$ which contains the
quasi--$\widetilde{\star}$--prime ideal $N'\cap D$.

\bf (2) \rm If $\star = {\star_{_{\!f}}}$ is e.a.b., then each $F \in
\boldsymbol{f}(D)$ is $\star$--e.a.b., thus
every  quasilocal   $\star$--overring (in particular, a
strong--$\star$--monolocality) is
necessarily a valuation domain, hence $\calL'(D, \star) =\calV(D,
\star)$.
Therefore, in this situation, $ \SKN(D, \star)=\Kr(D, \star)$.  Using
the
previous argument (and \cite[Proposition 4.1   (2)]{FL3}),   for each
semistar operation $\star$, passing to the
e.a.b. semistar operation of finite type $\star_{a}$, we have:
$$
  \SKN(D,  \star_{a}) =\Kr(D, \star_{a}) = \Kr(D, \star)\,.$$
\vskip -15 pt \hfill $\Box$     \vskip 12pt

%COROLLARY 5.5
\begin{coor} \label{co:1.16} \sl  Let $\star$ be a semistar operation
on an
integral domain $D$. Then $D$ is a P$\star$MD if and only if
$ \SKN(D,  \widetilde{\star})=\SKN(D,  \star_{a})$.
\end{coor}

\bf Proof.  \rm This statement is a straightforward consequence of
Proposition \ref{pr:1.12} and \cite[Theorem 3.1 and Remark
3.1]{FJS}.   \hfill $\Box$

\bigskip

We have defined and done some analysis on the generalized
Kronecker--Nagata rings using the (strong) monolocalities.  This was
motivated by characterizations of the classical Kronecker and Nagata
rings.  Both the Kronecker and Nagata rings have definitions involving
rational functions, content ideals, and semistar operations.  We turn
toward generalizing along these lines now.

%PROPOSITION 5.6
\begin{prro} \label{prop:1.20} \sl Let $\star$ be a semistar
operation on an
integral domain $D$ with quotient field $K$ and  let $\calL'
= \calL' (D, \star) $ be the set of all the
strong--$\star$--monolocalities of $D$.
Set:
$$\begin{array}{rl}
{\SKN}_{\#}(D, \star) \!:=
   \{ z\in\!  K(X) \!\mid & \hskip -5pt \forall \  L'
   \in\!  \calL',\ \exists\ g_{ L'}\in \!D[X],\ g_{ L'} \neq 0,
\mbox{ with } \\         & \hskip -3pt  zg_{ L'} \in
   \!D[X],\, \
    \boldsymbol{c}(zg_{ L'}) \subseteq
   \boldsymbol{c}(g_{ L'})L' \, \mbox{  and } \\ & \hskip -3pt
\boldsymbol{c}(g_{ L'})L'  \mbox{ is principal in }
   L^\prime\}\,.
     \end{array} $$
 Then
$\SKN(D, \star) \! = \SKN_{\#}(D, \star)$\,.
\end{prro}

\bf Proof.   \rm       Let $ z  \in
\SKN_{\#}(D,\star) $,  let $L' \in \calL'$ and let $g_{
L'}\in D[X]$ be such that $\boldsymbol{c}(g_{ L'})L'$ is a nonzero
principal ideal of $L'$
and $c(zg_{ L'}) \subseteq c(g_{ L'})L'$.  Set $f_{
L'}:= zg_{ L'}$.   Write $g_{ L'} :=a_{0}+ a_{1}X+\ \ldots\ +
a_{n}X^n \in D[X]$.  Since $L'$ is a
(strong--)$\star$--monolocality of $D$ and
$\boldsymbol{c}(g_{ L'})L'$ is principal then,  by the Claim in
the proof of
Lemma \ref{lm:1.4} (4), we have $
\boldsymbol{c}(g_{ L'})L'=  a_iL'$
for some $a_i$.
Hence, $\frac{f_{ L'}}{a_i} \in L'[X]$
and $\frac{g_{ L'}}{a_i} \in L'[X]$.  Moreover, $\frac{g_{ L'}}{a_i}$
is a primitive
polynomial  of $L'[X]$  (since one of its coefficients is a unit in
$L'$), hence:
\[
 z =  \frac{f_{ L'}}{g_{ L'}} =
\frac{\frac{f_{ L'}}{a_i}}{\frac{g_{ L'}}{a_i}} \in L'(X)\,.
\] Therefore, we have proven that $\SKN_{\#}(D, \star) \subseteq
\SKN(D, \star)$\,.

 In order to complete the proof, we need to show that
$\SKN(D, \star)
\subseteq \SKN_{\#}(D, \star) $.  If $z \in \SKN(D, \star)$,  then,
for each strong--$\star$--monolocality $L'$
of $D$, there exist $\varphi_{_{\!L'}},\ \psi_{_{\!L'}}\in L'[X]$,
with
$\psi_{_{\!L'}}\neq 0$, such that $ z =
\frac{\varphi_{_{\!L'}}}{\psi_{_{\!L'}}}$ and
$\boldsymbol{c}({\psi_{_{\!L'}}}) =L'$.  Therefore, we can find
$f_{_{\!L'}},\ g_{_{\!L'}} \in D[X]$ and two nonzero elements
$\alpha_{_{\!L'}},\
\beta_{_{\!L'}} \in D$ such that $f_{_{\!L'}} =
\alpha_{_{\!L'}}\varphi_{_{\!L'}}$, \ $
g_{_{\!L'}}=\beta_{_{\!L'}}{\psi_{_{\!L'}}}$ and thus:
$$  z =
\frac{ \varphi_{_{\!L'}} }
      {  \psi_{_{\!L'}} } =
\frac{
         \frac{ f_{_{\!L'}} }{ \alpha_{_{\!L'}} }
       }
      {
        \frac{ g_{_{\!L'}} }{ \beta_{_{\!L'}} }
      } =
  \frac{ \beta_{_{\!L'}} f_{_{\!L'}} } { \alpha_{_{\!L'}}
g_{_{\!L'}} }
      $$
      with
      $\boldsymbol{c}({f_{_{\!L'}}})L' \subseteq \alpha_{_{\!L'}} L'$
and
      $\boldsymbol{c}({g_{_{\!L'}}})L' = \beta_{_{\!L'}} L'$.
Therefore,
      $z \in  \SKN_{\#}(D,\star) $,    since
      $\boldsymbol{c}(\beta_{_{\!L'}}{f_{_{\!L'}}})L' =
      \beta_{_{\!L'}}\boldsymbol{c}({f_{_{\!L'}}})L' \subseteq
      \beta_{_{\!L'}}\alpha_{_{\!L'}}L '=
      \alpha_{_{\!L'}}\boldsymbol{c}({g_{_{\!L'}}})L' =
      \boldsymbol{c}(\alpha_{_{\!L'}}{g_{_{\!L'}}})L'$ and
      $\boldsymbol{c}(\alpha_{_{\!L'}}{g_{_{\!L'}}})L'$ is principal
in $L'$, for
      each $L \in \calL'$.  \hfill $\Box$

      \bigskip

Note that it follows immediately from the definition that
$D[X] \subseteq \SKN_{\#}(D,\star) \subseteq K(X)$.  Hence the 
quotient
field of $\SKN_{\#}(D,\star) =  \SKN(D,\star)$ is $K(X)$.

Our next result gives a basic property of $\SKN(D,\star)$ reminiscent
of Kronecker function ring and Nagata ring properties.

%PROPOSITION 5.7
\begin{prro}  \label{prop:5.8} \sl  Let $\star$ be a semistar
operation on an
integral domain $D$.
For each $J:= (a_{0},\ a_{1},\ \ldots,\ a_{n})D
\in \boldsymbol{f}(D)$, with $J\subseteq D$ and $J$ $\star$--e.a.b.,
let $g:=a_{0}+ a_{1}X+\ \ldots\ + a_{n}X^n \in D[X]$, then:

$$J\SKN(D, \star) =J^\star \SKN(D, \star)=g\SKN(D, \star)\,.$$
\end{prro}

{\bf Proof:}
 \rm  First note that, by definition, $J
=\boldsymbol{c}(g)D$.  Moreover,
for each $k$, with $0\leq k \leq n$, we have $a_k/g \in \SKN_{\#}(D,
\star)  \ (= \SKN(D,\star))$,  since
$\boldsymbol{c}(g)L' $ is principal in $L'$, for each $ L' \in
\calL'$.   Hence $J\SKN(D,\star))\subseteq g\SKN(D, \star)$.  Clearly,
$g \in J\SKN(D, \star)= \boldsymbol{c}(g)\SKN(D, \star)$.  It follows
that $J\SKN(D, \star) = g\SKN(D, \star)$.

On the other hand, let $\alpha =d/d' \in J^\star$, with $d, d' \in D$
and $d' \neq 0$.  Then, by Lemma
\ref{lm:1.4} (4), for each (strong--)$\star$--monolocality $L'$ of
$D$ we have  $dD
=d'\alpha D \subseteq d'\alpha L'
\subseteq d' J^\star L' = d' JL' = d' \boldsymbol{c}(g)L'$, thus
$\alpha/g = d'\alpha/d'g \in \SKN_{\#}(D,\star) \ (=
\SKN(D,\star))$,  since
$\boldsymbol{c}(g)L' =J\!L'$ is principal in $L'$,  for each $L' \in
\calL'$.
Hence, we have that  $J^\star \subseteq
g\SKN(D,\star)=J\SKN(D,\star)$, thus
$J^\star \SKN(D,\star) = J\SKN(D,\star)$. ~\hfill $\Box$

\bigskip

The rational function definition of the strong Kronecker--Nagata ring
is somewhat cumbersome.  We introduce now the notion of an   
``almost
e.a.b.--ideal''   in an effort to make the definition cleaner.

 % DEFINITION 5.8
\begin{deef} \label{rk:1.20} \rm

 Let $\star$ be a semistar operation on an
integral domain $D$  and let $F \in
\boldsymbol{f}(D)$, we say that $F$ is \it an
almost--$\star$--e.a.b.--ideal \rm if, for each 
strong--$\star$--monolocality
$L'$ of $D$,
  $F\!L'$ is
  a principal fractional ideal of $L'$.
\end{deef}

We collect in the following statement some of the basic properties of 
the almost--$\star$--e.a.b. ideals.

%PROPOSITION 5.9 
\begin{prro}   \sl  Let $\star$ be a semistar operation on an
integral domain $D$ and let $F \in \boldsymbol{f}(D) $   
\bara
\bf \item \sl
       If $F$ is $\star$--e.a.b. then $F$ is
almost--$\star$--e.a.b..

       \bf \item \sl   $F$ is  almost--$\star$--e.a.b.  if and only
if $F$
       is
       almost--$\star{_{_{\! f}}}$--e.a.b..

         \bf \item \sl If $F$ is  almost--$\star$--e.a.b. then
         $F\!L'$ is  $\star_{\iota}$--e.a.b.,  for each
strong--$\star$--monolocality
         $L'$ of $D$, with $\iota\ (= \iota_{_{\!  L'}}): D
\hookrightarrow L'$ being the canonical
embedding.

\bf \item \sl  Let $F := (a_{1},\ a_{2},\
    \ldots, a_{n})\!  D$ be an
    almost--${\star}$--e.a.b.--ideal,  then,  for each
strong--$\star$--monolocality
         $L'$ of $D$,
             $$
             F\! L'=F^\star\!  L' = (F\!  L')^{\star_{_{\!f}}} =
a_{i}L'\,,
             $$
             \ for some $i$, with $1\leq i \leq n$.
             \eara \rm
             \end{prro}

        \bf (1) \rm and \bf (2) \rm are obvious since if $F$ is
        $\star$--e.a.b., then $F\!L'$ is
        principal, for each strong--$\star$--monolocality $L'$ of
$D$, and since the
        strong--$\star$--monolocalities coincide with the
strong--$\star{_{_{\! f}}}$--monolocalities.

        \bf (3) \rm If $G,\ H \in \boldsymbol{f}(L')$ and $F \in
\boldsymbol{f}(D)$ is an
        almost--$\star$--e.a.b., then $F\!L' =zL'$,  for some
nonzero
        element $z$, and thus
        $(FG)^{\star_{\iota}} =(FL'G)^{\star}
        \subseteq  (F\!L'H)^{\star} = (FH)^{\star_{\iota}}$, then
       $ (zL'G)^\star
       \subseteq  (zL'H)^{\star}$.  Hence,
       $G^{\star_{\iota}} =G^\star  \subseteq
       H^\star =H^{\star_{\iota}} $.

       \bf (4) \rm  This statement
is a consequence of the Claim in  the proof
       of Lemma \ref{lm:1.4}  (4).  \hfill $\Box$

\bigskip

This allows us to state a new definition.

%DEFINITION 5.10
\begin{deef} \label{def:5.10} \rm     Let $\star$
 be a semistar operation on an
             integral domain $D$   with quotient field $K$.     
Then we define
  $$\begin{array}{rl}
\SKN_{\boldsymbol{c}}(D, \star) :=
    \{ z\in\!  K(X) \!\mid & \hskip -7pt \exists \ g \in \!D[X], \ g
\neq 0, \mbox{ with } zg \in
   \!D[X],\   \boldsymbol{c}(zg) \subseteq
   \boldsymbol{c}(g)^{\star}
    \mbox{  and } \\ 
    &  \hskip -6pt 
   \boldsymbol{c}(g) \mbox{  is an almost--$\star$--e.a.b. ideal 
of  }
   D \} \,,
     \end{array}$$
     $$
\begin{array}{rl}
    \KN_{\boldsymbol{c}}(D, \star) :=  \{ z\in K(X) \mid & \exists
\  g \in
D[X],\
  g \neq 0,\ \mbox{ such that } zg \in D[X]\,, \ \mbox{ and } \\
   &   \boldsymbol{c}(zg) \subseteq \boldsymbol{c}(g)^\star,\
   \boldsymbol{c}(g)\ \mbox{ is a }\star\mbox{--e.a.b. ideal} \}\,.
     \end{array}
     $$

     \end{deef}

  It is clear that, in general, we have
  $\SKN_{\boldsymbol{c}}(D, \star) \subseteq \SKN(D, \star) \ (= 
\SKN_{\#}(D, \star)).$

     In fact, note that,   for each $g \in \!D[X],\ g \neq 0,$ and
      for each $L'\in \calL'$,
     $ \boldsymbol{c}(g)^{\star}L' = \boldsymbol{c}(g)L'$;
     moreover,
     $\boldsymbol{c}(g)$ is an   almost--$\star$--e.a.b. ideal of
     $D$ if and only if (by definition)
     it is a principal ideal of $L'$, for each $L'\in \calL'$.
     \smallskip

 We suspect that in fact
 $\SKN_{\boldsymbol{c}}(D, \star) = \SKN(D, \star)$, but we do not 
have
 a proof.  We do demonstrate below that 
 $\KN_{\boldsymbol{c}}(D, \star) = \KN(D, \star)$.

 \bigskip

 We turn now to investigating the properties of $\KN(D,\star)$.  With
 this ring    we will have   more luck demonstrating properties that
 reflect the classical properties of the Kronecker function rings and
 Nagata rings.  In particular, when we localize a Kronecker
 function ring $\Kr(D,\star)$ at a maximal ideal we obtain $V(X)$ for
 some  ($\star$--)valuation   overring $V$ of $D$.  Similarly, 
when we localize
 a Nagata ring  $\Na(D,\star)$ at a maximal ideal we obtain $D_{Q}(X)$
 for some  (quasi--$\star$--)prime    ideal $Q$ of $D$.  We 
obtain similar results with
 $\KN(D,\star)$.

 %THEOREM 5.11
\begin{thee} \label{th:1.21} \sl Let $\star$
   be a semistar operation on an
integral domain $D$ with quotient field $K$.  \bara
%1
\bf \item \sl 
   $ \KN_{\boldsymbol{c}}(D, \star) $
     is an integral domain with quotient field $K(X)$.
     %2
     \bf \item \sl $ \Na(D, \star) \subseteq  \KN_{\boldsymbol{c}}(D,
\star)  \subseteq \SKN_{\boldsymbol{c}} (D, \star)
     \ (\subseteq  \SKN_{\#}(D, \star) =\SKN(D, \star))\,.$
     %3
     \bf \item \sl For each $J:= (a_{0},\ a_{1},\ \ldots,\  a_{n})D
\in
     \boldsymbol{f}(D)$, with $J\subseteq D$ and $J$
$\star$--e.a.b.,  let
     $g:=a_{0}+ a_{1}X+\ \ldots\ + a_{n}X^n \in D[X]$, then:
     $$
     J\KN_{\boldsymbol{c}}(D, \star) =J^\star\KN_{\boldsymbol{c}}(D,
\star)=g\KN_{\boldsymbol{c}}(D, \star)\,.$$
     %4
     \bf \item \sl For each  prime   ideal ${
\boldsymbol{\mathfrak{p}} }$ of\  $\KN_{\boldsymbol{c}}(D, \star)$
     and for each
     $J:= (a_{0},\ a_{1},\ \ldots,\  a_{n})D \in
     \boldsymbol{f}(D)$, with $J\subseteq D$ and $J$
$\star$--e.a.b.,  let
     $g:=a_{0}+ a_{1}X+\ \ldots\ + a_{n}X^n \in D[X]$, then there
exists an
     index $i$, with $0 \leq i \leq n$, such that:
     $$
     J\KN_{\boldsymbol{c}}(D, \star)_{{ \boldsymbol{\mathfrak{p}}
}} =
     J^\star\KN_{\boldsymbol{c}}(D, \star)_{{
\boldsymbol{\mathfrak{p}} }}
     =g\KN_{\boldsymbol{c}}(D, \star)_{{ \boldsymbol{\mathfrak{p}}
}} =
     a_{i}\KN_{\boldsymbol{c}}(D, \star)_{{
\boldsymbol{\mathfrak{p}} }}\,.$$

       \hskip -28pt For each  prime  ideal ${
\boldsymbol{\mathfrak{p}} }$ of $\KN_{\boldsymbol{c}}(D, \star)$,
set
       $L({ \boldsymbol{\mathfrak{p}} }) :=
       \KN_{\boldsymbol{c}}(D, \star)_{{ \boldsymbol{\mathfrak{p}}
}} \cap K$.

       %5
        \bf \item \sl For each  prime  ideal ${
\boldsymbol{\mathfrak{p}} }$ of $\KN_{\boldsymbol{c}}(D, \star)$,
then
     $L({ \boldsymbol{\mathfrak{p}} })$ is a
$\star$--monolocality of
     $D$ (with maximal ideal\  
$\boldsymbol{\boldsymbol{\mathscr{P}}}  := {
\boldsymbol{\mathfrak{p}} }\!  \KN_{\boldsymbol{c}}(D,
     \star)_{{ \boldsymbol{\mathfrak{p}} }} \cap
     L({ \boldsymbol{\mathfrak{p}} })$).

     %6
     \bf \item \sl For each  prime  ideal ${
\boldsymbol{\mathfrak{p}} }$ of $\KN_{\boldsymbol{c}}(D, \star)$,
    the localization $\KN_{\boldsymbol{c}}(D,
     \star)_{{ \boldsymbol{\mathfrak{p}} }}$
     coincides with the Nagata ring $L({ \boldsymbol{\mathfrak{p}}
}) (X)$ (with maximal ideal $ {\boldsymbol{\mathscr{P}}}(X):= {\boldsymbol{\mathscr{P}}}L({ \boldsymbol{\mathfrak{p}}}) (X)$) 
     and ${ \boldsymbol{\mathfrak{p}} }$ coincides with 
    $ {\boldsymbol{\mathscr{P}}}(X)  \cap \KN_{\boldsymbol{c}}(D, 
\star)$.

     %7
     \bf \item \sl
   Every $\star$--monolocality of an
    integral domain $D$ contains a minimal $\star$--mono\-locality of
    $D$. If we denote by ${\calL}(D, \star)_{min}$, or simply
    by
       ${\calL}_{min}$,  the set of all the minimal
       $\star$--monolocalities of $D$, then
       ${\calL}(D, \star)_{min} =
       \{ L({\boldsymbol{\mathfrak{m}}}) \mid
       {\boldsymbol{\mathfrak{m}}} \in \Max(\KN_{\boldsymbol{c}}(D,
\star)) \}$ and
       $$ \KN_{\boldsymbol{c}}(D, \star) ={\KN}(D, \star) = \cap
\{L(X) \mid L
    \in {\calL}_{min} \}\,.$$
     In particular, $\Na(D, \star)
\subseteq {\KN}(D, \star) \ (\subseteq {\SKN}(D, \star) \subseteq
\Kr(D, \star))$\,.

\eara
    \end{thee}

     \bf Proof.  \rm
     \bf (1) \rm    Note that:

     \bf Claim 1. \sl If $g$ and $ h$ are two nonzero
     polynomials of  $D[X]$ and $\boldsymbol{c}(g)$ is a
$\star$--e.a.b.
     ideal of $D$ then
 $(\boldsymbol{c}(g) \boldsymbol{c}(h))^\star=
 \boldsymbol{c}(gh)^\star$.  Furthermore, if
 $\boldsymbol{c}(h)$ is also a $\star$--e.a.b. ideal of $D$, then
 $\boldsymbol{c}(gh)$ is a $\star$--e.a.b. ideal of
 $D$. \rm

 The previous claim is straightforward consequence of the
Dedekind-Mertens
 Lemma \cite[Theorem 28.1]{G} and of the definition of
$\star$--e.a.b. ideal.

     Let $z:= f/g, \ z':= f'/g' \in \KN_{\boldsymbol{c}}(D, \star) $,
with
     $\boldsymbol{c}(f) \subseteq  \boldsymbol{c}(g)^\star$,
     $\boldsymbol{c}(f') \subseteq  \boldsymbol{c}(g')^\star$ and
     $\boldsymbol{c}(g)$ and $\boldsymbol{c}(g')$
     $\star$--e.a.b. ideals of $D$.  From Claim 1, we deduce
immediately that $zz'=ff'/gg' \in  \KN_{\boldsymbol{c}}(D, \star) 
$.   

  In order to see that $z-z'$ belongs to $\KN_{\boldsymbol{c}}(D,
\star) $, it is sufficient to observe that $z-z' = (fg'
   -f'g)/gg'$ and
  $$\boldsymbol{c}(fg'-f'g) \subseteq \boldsymbol{c}(fg'-f'g)^\star
  \subseteq (\boldsymbol{c}(fg')^\star +\
\boldsymbol{c}(f'g)^\star)^\star
  \subseteq \boldsymbol{c}(gg')^\star\,.$$

  Clearly, \ $D[X] \subseteq \KN_{\boldsymbol{c}}(D, \star)\subseteq
K(X)$, hence the
quotient field of $\KN_{\boldsymbol{c}}(D, \star) $ is $K(X)$.

     %2
 \bf (2) \rm   To prove that  $\Na(D, \star) \subseteq
\KN_{\boldsymbol{c}}(D, \star)$, note that the definition of 
$\KN_{\boldsymbol{c}}(D, \star)$ generalizes the definition of 
 $\Na(D, \star)$ (Definition \ref{def:3.1}) by replacing 
$\star_f$--invertible 
 ideals with the larger class of $\star$--e.a..b. ideals.  The result 
is then clear.
The second inclusion  is an
easy consequence of the fact that, if $\boldsymbol{c}(zg) \subseteq
 \boldsymbol{c}(g)^\star$ and $
   \boldsymbol{c}(g)\ \mbox{ is  }\star\mbox{--e.a.b. }$ then, for
each $L \in\calL$,
    $\boldsymbol{c}(g)^\star L = \boldsymbol{c}(g)L$ is a principal
    ideal of $L$ (Lemma \ref{lm:1.4}  (4)).     Therefore
    we have $ \KN_{\boldsymbol{c}}(D, \star) \subseteq
\SKN_{\boldsymbol{c}}(D, \star)
  $\  (Definition
     \ref{rk:1.20}).

  %3
\bf   (3) \rm  
      \it Mutatis mutandis \rm the proof ot the   equality    $
J\KN_{\boldsymbol{c}}(D, \star)$ = $
\boldsymbol{c}(g)\KN_{\boldsymbol{c}}(D, \star) $ = $
g\KN_{\boldsymbol{c}}(D, \star)$
 is analogous to the
proof  of Proposition
\ref{prop:5.8}.   

      More precisely, first, note that by
definition $J =\boldsymbol{c}(g)D$.
   Moreover,  for each $k$, with $0\leq k \leq n$, we have $a_k/g \in
\KN_{\boldsymbol{c}}(D, \star)$.  Hence $J\KN_{\boldsymbol{c}}(D,
\star) \subseteq g\KN_{\boldsymbol{c}}(D, \star)$.  Clearly, $g
   \in \boldsymbol{c}(g)\KN_{\boldsymbol{c}}(D,
\star)=J\KN_{\boldsymbol{c}}(D, \star)$,   and so 
  $J\KN_{\boldsymbol{c}}(D, \star) $ = $
g\KN_{\boldsymbol{c}}(D, \star)$.

 On the other hand, let $\alpha:= d/d' \in J^\star$, with $d, d' \in
 D$, $d' \neq 0$ Then $\alpha/g = d'\alpha/d'g = d/d'g \in
\KN_{\boldsymbol{c}}(D, \star)$, since $dD \subseteq d'J^\star
=d'\boldsymbol{c}(g)^\star$.  Hence,
$J^\star \subseteq g\KN_{\boldsymbol{c}}(D, \star) $ = $
J\KN_{\boldsymbol{c}}(D, \star)$ and so $J\KN_{\boldsymbol{c}}(D,
\star)={J^\star}\KN_{\boldsymbol{c}}(D, \star)$.

      %4
\bf (4) \rm All the equalities follow trivially from (3) except the
last one,
 involving $a_i$.  This equality holds  because
 $\KN_{\boldsymbol{c}}(D, \star)_{ \boldsymbol{\mathfrak{p}} }$
is  quasilocal and  
 $g\KN_{\boldsymbol{c}}(D, \star)_{ \boldsymbol{\mathfrak{p}}}= 
J\KN_{\boldsymbol{c}}(D, \star)_{ \boldsymbol{\mathfrak{p}}} $ $= (a_0, a_1, \ldots, a_n)\KN_{\boldsymbol{c}}(D, \star)_{ \boldsymbol{\mathfrak{p}}}$
 is a principal.  By a standard technique
  (Claim in
 the proof of Lemma \ref{lm:1.4}  (4)),   an invertible ideal of 
a  quasilocal   
domain
  which is generated by a finite list of elements is actually
generated by
 one of those elements.

      %5
 \bf (5) \rm  It is clear that $L({ \boldsymbol{\mathfrak{p}}
})$ is a
 quasilocal    overring of $D$, with maximal ideal   
$\boldsymbol{\mathscr{P}}$.
In
 order to show that $L({ \boldsymbol{\mathfrak{p}} })$ is a
$\star$--monolocality of
     $D$,
take $J:= (a_{0},\ a_{1},\ \ldots,\ a_{n})D $ which is a nonzero
$\star$--e.a.b. ideal of $D$.  It is clear that
$a_k L({ \boldsymbol{\mathfrak{p}} }) \subseteq
J^\star L({ \boldsymbol{\mathfrak{p}} })  $, for each $0 \leq k
\leq n$.  Let $\alpha \in J^\star$.  Since, by (4),
 $J^\star \KN_{\boldsymbol{c}}(D, \star)_{
\boldsymbol{\mathfrak{p}} } =
 a_i \KN_{\boldsymbol{c}}(D, \star)_{ \boldsymbol{\mathfrak{p}}
}$,
for some $i$, then
$\alpha/a_i \in \KN_{\boldsymbol{c}}(D, \star)_{
\boldsymbol{\mathfrak{p}} } \cap K =
L({ \boldsymbol{\mathfrak{p}} })$.
Therefore,
 $J^\star \subseteq a_i L({ \boldsymbol{\mathfrak{p}} })
\subseteq JL({ \boldsymbol{\mathfrak{p}} }) \subseteq
J^\star L({ \boldsymbol{\mathfrak{p}} })$ and so
$ a_iL({ \boldsymbol{\mathfrak{p}} })= J^\star L({
\boldsymbol{\mathfrak{p}} })$.

      %6
\bf (6) \rm  We start by proving the following:

\bf Claim 2. \sl  Let $L$ be a $\star$--monolocality of $D$,
then  $\KN_{\boldsymbol{c}}(D, \star)\subseteq L(X)$.  In
particular,  $\KN_{\boldsymbol{c}}(D, \star) \subseteq
{\KN}(D, \star)$. \rm

Let $f/g \in \KN_{\boldsymbol{c}}(D, \star)$ with $\boldsymbol{c}(g)$
a $\star$--e.a.b.
ideal of $D$
and $\boldsymbol{c}(f) \subseteq \boldsymbol{c}(g)^\star$.  Write $g :=a_{0}+ a_{1}X+\ \ldots\ +
a_{n}X^n \in D[X]$.  Since $L$ is a $\star$--monolocality of $D$
then, by
Lemma \ref{lm:1.4} (4), we have:
\[
\boldsymbol{c}(g)L= \boldsymbol{c}(g)^\star L =
(\boldsymbol{c}(g)L)^{\star_{_{\!f}}} = a_iL
\]
for some $a_i$.
Hence, $\frac{f}{a_i} \in L[X]$
and $\frac{g}{a_i} \in L[X]$.  Moreover, $\frac{g}{a_i}$ is a
primitive
polynomial   in $L[X]$  (since one of its coefficients is a unit in
$L$).  Hence
\[
\frac{f}{g} = \frac{\frac{f}{a_i}}{\frac{g}{a_i}} \in L(X)\,,
\]
 and so Claim 2 is proven.

Note that, by (5) and Claim 2, we have
 $L({ \boldsymbol{\mathfrak{p}} })(X)
\supseteq \KN_{\boldsymbol{c}}(D, \star)$.  Note also that, since $
L({ \boldsymbol{\mathfrak{p}} }) \subseteq
\KN_{\boldsymbol{c}}(D, \star)_{ \boldsymbol{\mathfrak{p}} }$
and $X \in \KN_{\boldsymbol{c}}(D, \star)$,
then $L({ \boldsymbol{\mathfrak{p}} })[X] \subseteq
\KN_{\boldsymbol{c}}(D, \star)_{ \boldsymbol{\mathfrak{p}} }$
and hence
$ \boldsymbol{\mathscr{P}}[X]   =
\left( \boldsymbol{\mathfrak{p}} \KN_{\boldsymbol{c}}(D,
\star)_{ \boldsymbol{\mathfrak{p}} }
\cap L({ \boldsymbol{\mathfrak{p}} })\right)[X]
=  \boldsymbol{\mathfrak{p}} \KN_{\boldsymbol{c}}(D,
\star)_{ \boldsymbol{\mathfrak{p}} }
\cap \left(L({ \boldsymbol{\mathfrak{p}} })[X]\right)$,
recalling that
$ \boldsymbol{\mathscr{P}} =
 \boldsymbol{\mathfrak{p}} \KN_{\boldsymbol{c}}(D, \star)_{
\boldsymbol{\mathfrak{p}} }
\cap L({ \boldsymbol{\mathfrak{p}} })$.
Clearly, $
\boldsymbol{\mathscr{P}}(X)  \cap \KN_{\boldsymbol{c}}(D, \star)$ is a
proper prime ideal of $\KN_{\boldsymbol{c}}(D, \star)$.

\bf Claim 3. \sl With the notation introduced above, $
\boldsymbol{\mathscr{P}}(X)  \cap \KN_{\boldsymbol{c}}(D, \star) 
\subseteq
 \boldsymbol{\mathfrak{p}} $. \rm

Let $\varphi \in  \boldsymbol{\mathscr{P}}(X)  \cap
\KN_{\boldsymbol{c}}(D, \star)$.  Then, we can write $\varphi= h /k$
where $h, k \in
L({ \boldsymbol{\mathfrak{p}} })[X]$ and
$h \in  \boldsymbol{\mathscr{P}}[X]$  and $k$ is primitive in
$L({ \boldsymbol{\mathfrak{p}} })[X]$.  We can also
write $\varphi = f/g$ where $f, g \in D[X]$ , $g \neq 0$, 
$\boldsymbol{c}(f) \subseteq
\boldsymbol{c}(g)^\star$ and
$\boldsymbol{c}(g)$ is a $\star$--e.a.b. ideal of $D$.
Since $L({ \boldsymbol{\mathfrak{p}} })$ is a
$\star$--monolocality
of $D$ (by (5)), it then follows
from Lemma \ref{lm:1.4} (4) that $g$ has a coefficient $a_i$ such that
$\boldsymbol{c}(g)L({ \boldsymbol{\mathfrak{p}} }) = a_i L({
\boldsymbol{\mathfrak{p}} })$ and, hence,
$\boldsymbol{c}(f) \subseteq a_iL({ \boldsymbol{\mathfrak{p}}
})$.  Then:
\[
\frac{f}{g} = \frac {\frac{f}{a_i}}{\frac{g}{a_i} } = \frac{h}{k}
\]
with $\frac{f}{a_i}, \frac{g}{a_i} \in L({
\boldsymbol{\mathfrak{p}} })[X]$.
Therefore
\[
k \frac{f }{a_i} = h \frac{g}{a_i}\,.
\]
Since $k $ and $\frac{g}{a_i} $ are primitive in
$L({ \boldsymbol{\mathfrak{p}} })[X]$ and $h \in
 \boldsymbol{\mathscr{P}}[X] $, then we must have $\frac{f}{a_i}\in 
\boldsymbol{\mathscr{P}}[X]  \subseteq
{ \boldsymbol{\mathfrak{p}} }\KN_{\boldsymbol{c}}(D,
\star)_{ \boldsymbol{\mathfrak{p}} }$.
Since $\frac{g}{a_i} $ is a unit in
$\KN_{\boldsymbol{c}}(D, \star)_{ \boldsymbol{\mathfrak{p}} }$
this implies that $f/g\ (=h/k = \varphi)$ belongs to
$
{ \boldsymbol{\mathfrak{p}} }\KN_{\boldsymbol{c}}(D,
\star)_{ \boldsymbol{\mathfrak{p}} } \cap
\KN_{\boldsymbol{c}}(D, \star) =
{ \boldsymbol{\mathfrak{p}} }$.

 We conclude the proof of (6).  By Claim 3,
the prime ideal $ {\boldsymbol{\mathfrak{p'}}}  :=
\boldsymbol{\mathscr{P}}(X) \cap
\KN_{\boldsymbol{c}}(D, \star)$ is
contained in  ${ \boldsymbol{\mathfrak{p}} }$,  thus   it is
clear that:
$$ L({ \boldsymbol{\mathfrak{p}} })[X]
\subseteq
\KN_{\boldsymbol{c}}(D, \star)_{ \boldsymbol{\mathfrak{p}} } \
\subseteq
\   \KN_{\boldsymbol{c}}(D, \star)_{
\boldsymbol{\mathfrak{p'}}} \ \subseteq \
L({ \boldsymbol{\mathfrak{p}} })(X) =
L({ \boldsymbol{\mathfrak{p}} })[X]_{ 
\boldsymbol{\mathscr{P}}[X]}\,.$$
Moreover, \
$ \boldsymbol{\mathfrak{p}} \KN_{\boldsymbol{c}}(D, \star)_{
\boldsymbol{\mathfrak{p}} }
\cap \left(L({ \boldsymbol{\mathfrak{p}} })[X]\right) = 
\boldsymbol{\mathscr{P}}[X]$.
Therefore,
$$ L({ \boldsymbol{\mathfrak{p}} })[X]_{ \boldsymbol{\mathscr{P}}[X]}
=
\KN_{\boldsymbol{c}}(D, \star)_{ \boldsymbol{\mathfrak{p}} } =
\KN_{\boldsymbol{c}}(D, \star)_{ \boldsymbol{\mathfrak{p'}}}
=
L({ \boldsymbol{\mathfrak{p}} })(X)\,,$$
hence we deduce that
${ \boldsymbol{\mathfrak{p'}}}={ \boldsymbol{\mathfrak{p}}
}$.

  %7
\bf (7) \rm is an easy consequence of (5) and (6). In fact, let $L
\in {\calL}(D, \star)$ and let ${N}$ be the maximal ideal of
$L$. We know by Claim 2 that $L(X) \supseteq \KN_{\boldsymbol{c}}(D,
\star)$. Set
${\boldsymbol{\mathfrak{n}}}: ={N}(X) \cap  \KN_{\boldsymbol{c}}(D,
\star)$. Let
${\boldsymbol{\mathfrak{m}}}$ be a maximal ideal of
$\KN_{\boldsymbol{c}}(D, \star)$ containing the prime ideal
${\boldsymbol{\mathfrak{n}}}$.
Then, by (5) and (6), we know  that:
$$
L({\boldsymbol{\mathfrak{m}}})(X)=\KN_{\boldsymbol{c}}(D,
\star)_{\boldsymbol{\mathfrak{m}}}
\subseteq  \KN_{\boldsymbol{c}}(D, \star)_{\boldsymbol{\mathfrak{n}}}
\subseteq
L(X)\,.$$
This fact implies that ${\calL}(D, \star)_{min} =
       \{ L({\boldsymbol{\mathfrak{m}}}) \mid
       {\boldsymbol{\mathfrak{m}}} \in \Max(\KN_{\boldsymbol{c}}(D,
\star)) \}$, and
        so:
       $$
       \begin{array}{rl}
          \KN_{\boldsymbol{c}}(D, \star) =& \hskip -4pt \cap \{
          \KN_{\boldsymbol{c}}(D, \star)_{\boldsymbol{\mathfrak{m}}}
\mid
           {\boldsymbol{\mathfrak{m}}} \in
\Max(\KN_{\boldsymbol{c}}(D, \star))\} =\\
           =& \hskip -4pt \cap \{
           L({\boldsymbol{\mathfrak{m}}}) (X) \mid
           {\boldsymbol{\mathfrak{m}}} \in 
\Max(\KN_{\boldsymbol{c}}(D,
\star))\} =\\
            =& \hskip -4pt \cap \{
           L(X) \mid
           L \in {\calL}(D, \star)_{min}\} =\\
            =& \hskip -4pt \cap \{
           L(X) \mid
           L \in {\calL}(D, \star)\} ={\KN}(D,
           \star)\,.
           \end{array}
           $$

      The last statement of (7) follows from (2).
        \hfill
$\Box$

%%%%%%%%%%%%% SECTION 6
\section{New semistar operations}

Given a semistar operation $\star$ on a domain $D$ we have associated
two collections of overrings $\calL  \ (=\calL(D, \star))$   and 
$\calL'  \ (={\calL}'(D, \star))$ and using these
collections we have constructed two rings of rational functions
$\KN(D,\star)$ and $\SKN(D,\star)$.  We can use these two collections
of overrings and two rings of rational functions to construct
four new semistar operations associated to $\star$.

%DEFINITION 6.1
\begin{deef}  \label{def:6.1} \rm  Let $D$ be a domain with quotient 
field $K$ and 
$\star$
a semistar operation
on $D$.  We define new semistar operations  on $D$    as 
follows.  For each $E \in  \overline{\boldsymbol{F}}(D)$,   
\begin{enumerate}
\item[\bf(a)] \it $ \wedge_{\calL'}$ \ defined by \   $E^{ 
\wedge_{\calL'}} =:
\cap \{EL' \mid L' \in \calL' \}$\,;
\item[\bf(b)] \it   $ \wedge_{\calL}$ \  defined by \  $E^{ 
\wedge_{\calL}} =:
\cap \{EL \mid L \in \calL \}$\';
\item[\bf(c)] \it  $\star_{_{\!\ell}}$ \ defined by \ 
$ E^{\star_{_{\!\ell}}}:= E\KN(D, \star) \cap K$\,;
\item[\bf(d)] \it  $\star_{_{\!\ell'}}$ \ defined by \ 
$ E^{\star_{{\!\ell'}}}:= E\SKN(D, \star) \cap K$\,.
\end{enumerate}

\end{deef}

Next we give some simple relations between these operations.

%PROPOSITION 6.2
\begin{prro} \label{pr:6.2}
Let $D$ be a domain and $\star$ a semistar operation on
$D$.  Then ${\star_{_{\!\ell}}}$ and $\star_{_{\!\ell'}}$ are semistar
operations of finite type of $D$
  and
  $$         \begin{array}{rl}     \wedge_{\calL} \leq& \hskip
-5pt  \wedge_{\calL'}\,,\;\;      {\star_{_{\!\ell}}} \leq
\wedge_{\calL} \,,\;\;  {\star_{_{\!\ell'}}} \leq  \wedge_{\calL'}
\\       & \hskip -5pt  {\star_{_{\!\ell}}}
\leq \star_{_{\!\ell'}} \,,\;\;   \star_{_{\!f}} \leq
\wedge_{\calL'}\,.       \end{array}
  $$
\end{prro}

{\bf Proof:} It is easy to verify that, for each integral domain $R$
with quotient field $K(X)$, such that $D \subseteq R\cap K$, the
operation  ${\circlearrowleft_R}$   defined by   
$E^{\circlearrowleft_R}   := ER \cap K$,  
for each $ E\in \overline{\boldsymbol{F}}(D)$ is a semistar operation
of finite type on $D$.  As a matter of fact, note that, for each
nonzero element $x\in K$ and for each $E \in
\overline{\boldsymbol{F}}(D)$, we have $(xER) \cap K = x(ER\cap K)$,
$ER = \bigcup \{ FR \mid F
    \subseteq E\,,\; F \in \boldsymbol{f}(D) \}$ and $( \bigcup \{ FR
\mid F
    \subseteq E\,,\; F \in \boldsymbol{f}(D) \}) \cap K = \bigcup \{ 
FR
\cap K  \mid F
    \subseteq E\,,\; F \in \boldsymbol{f}(D) \}$.   Therefore, in
particular, ${\star_{_{\!\ell}}}  \  (= 
\circlearrowleft_{{\footnotesize  \mbox{KN}}(D, \star)})$   and 
${\star_{_{\!\ell'}}}   \ (= \circlearrowleft_{{\footnotesize 
\mbox{KN}'}(D, \star)})$   are
semistar operations of finite type on $D$.
For each strong--$\star$--monolocality $L'$ of $D$, we have
 that $F^\star \subseteq F\!L'$, thus $F^{\star}\!L'=
 F\!L'$, hence in particular $F^{\wedge_{\calL'}}
=(F^\star)^{\wedge_{\calL'}}$.  Therefore, $ \star_{_{\!f}}  \leq
\wedge_{\calL'}$.  Moreover,
 $$
 \begin{array}{rl}
     F\!\SKN(D, \star)\cap K =& \hskip -5pt  \left(F(\bigcap \{L'(X)
\mid L' \in \calL'\})\right) \cap
 K  \\ \subseteq &\hskip -5pt \left(\bigcap \{F\!L'(X)  \mid L' \in
\calL'\}\right) \cap
 K\\ &\hskip -5pt \bigcap \{F\!L'(X) \cap
 K \mid L' \in \calL'\} =  \bigcap \{F\!L' \mid L' \in \calL'\}\,.
\end{array}
 $$  We conclude that ${\star_{_{\!\ell'}}} \leq  \wedge_{\calL'}$.
Similarly, it can be shown that
${\star_{_{\!\ell}}} \leq
\wedge_{\calL}$.  Finally, since
$\ \KN(D, \star)  \subseteq \SKN(D, \star) $  (Proposition 
\ref{lm:1.7}), then
$ {\star_{_{\!\ell}}} \leq {\star_{_{\!\ell'}}}$.   \hfill $\Box$

\bigskip

If we restrict to just ${\star_{_{\!\ell}}}$ and $\wedge_{\calL}$ we
can prove more.

%PROPOSITION 6.3
\begin{prro} \label{pro:6.3}  Let $D$ be a domain and 
$\star$ a semistar operation on
$D$.
 Then  $ \star_{_{\!\ell}}$ of $D$ satisfies
$$ \widetilde{\star}  \leq {\star_{_{\!{\ell}}}} = \wedge_{\calL}
\leq  \ {\star}_{_{\!{f}}} \,.$$

\end{prro}

{\bf Proof:}
%8
      For each $E \in
\overline{\boldsymbol{F}}(D)$:      $$      \begin{array}{rl}
E^{\star_{_{\!\ell}}}= & \hskip -5 pt E\KN(D, \star)   \cap K =
\left(\cap \{E\KN(D, \star)_{\boldsymbol{\mathfrak{m}}}    \mid
{\boldsymbol{\mathfrak{m}}} \in \Max(\KN(D, \star) \}\right) \cap K
\\      = & \hskip -5 pt \cap \{E\KN(D,
\star)_{\boldsymbol{\mathfrak{m}}} \cap K \mid
{\boldsymbol{\mathfrak{m}}} \in \Max(\KN(D, \star) \} \\

= & \hskip -5 pt \cap \{EL({\boldsymbol{\mathfrak{m}}}) (X)\cap K
\mid {\boldsymbol{\mathfrak{m}}} \in \Max(\KN(D, \star) \}
\\            = & \hskip -5 pt  \cap
\{EL({\boldsymbol{\mathfrak{m}}}) \mid {\boldsymbol{\mathfrak{m}}}
\in \Max(\KN(D, \star) \} = E^{\wedge_{\calL_{min}}} =
E^{\wedge_{\calL}} \,.            \end{array}
$$

    Let $J:=(a_{1}, a_{2}, \ldots, a_{n})D \in
\boldsymbol{f}(D)$.   Suppose that
$\alpha \in J^{\star_{_{\!{\ell}}}}$.  Since $\KN_{\boldsymbol{c}}(D,
\star) = \KN(D, \star)$, then in $\KN(D, \star)$ we can write:
\[
\alpha = a_{1}\!\left(\frac{f_{1}}{g_{1}}\right) +
a_{2}\!\left(\frac{f_{2}}{g_{2}}\right) +
\cdots + a_{n}\!\left(\frac{f_{n}}{g_{n}}\right)
\]
with $\boldsymbol{c}(f_{i}) \subseteq \boldsymbol{c}(g_{i})^{\star}$
and $\boldsymbol{c}(g_{i})$ a $\star$--e.a.b. ideal of
$D$,  for each $i$, $1\leq i \leq n$.

Let $g :=  g_{1}g_{2} \cdots g_{n}$ and, for each $i$,  let
$\check{g}_{i} := g_{1}g_{2} \cdots g_{i -1}g_{i+1}
\cdots g_{n} = g/g_{i}$.  Then we can write:
\[
g\alpha = a_{1}f_{1}\check{g}_{1} + a_{2}f_{2}\check{g}_{2}  + \cdots
+ a_{n}f_{n}\check{g}_{n}\,.
\]

Therefore
 \[
 \begin{array}{rl}
    \boldsymbol{c}(g\alpha)^{\star} =& \hskip -4pt
\boldsymbol{c}(a_{1}f_{1}\check{g}_{1} + a_{2}f_{2}\check{g}_{2}  +
\cdots
+ a_{n}f_{n}\check{g}_{n})^{\star} \subseteq \\
\subseteq  &\hskip -4pt
\boldsymbol{c}(a_{1}f_{1}\check{g}_{1})^{\star} +
\boldsymbol{c}(a_{2}f_{2}\check{g}_{2})^{\star} +\cdots
+ \boldsymbol{c}(a_{n}f_{n}\check{g}_{n})^{\star} \subseteq  \\
 \subseteq  & \hskip -4pt \boldsymbol{c}(a_{1}g_{1}g_{2} \cdots
g_{n})^{\star} + \boldsymbol{c}(a_{2}g_{1}g_{2}
\cdots
g_{n})^{\star} +
\cdots + \boldsymbol{c}(a_{n}g_{1}g_{2} \cdots g_{n})^{\star}
\subseteq \\
 \subseteq  & \hskip -4pt  (a_{1}, a_{2}, \cdots,
a_{n})\boldsymbol{c}(g_{1}g_{2}\cdots
g_{n})^{\star} =
(a_{1}, a_{2}, \cdots, a_{n})\boldsymbol{c}(g)^{\star} \subseteq \\
 \subseteq  & \hskip -4pt
\left((a_{1}, a_{2}, \cdots, a_{n})\boldsymbol{c}(g)\right)^{\star}\,.
\end{array}
\]

Since each $\boldsymbol{c}(g_{i})$ is a $\star$--e.a.b. ideal of $D$,
then we know that
$\left(\boldsymbol{c}(g_{1})\boldsymbol{c}(g_{2})\cdots
\boldsymbol{c}(g_{n})\right)^{\star} = \boldsymbol{c}(g_{1}g_{2}
\cdots
g_{n})^{\star} = \boldsymbol{c}(g)^{\star}$ and that
$\boldsymbol{c}(g)$ is
a $\star$--e.a.b. ideal of $D$ (Claim 1  in the proof of Theorem 
\ref{th:1.21}).     It follows that:
\[
(\boldsymbol{c}(g)\alpha)^{\star}
\subseteq
\left((a_{1}, a_{2}, \cdots
, a_{n})\boldsymbol{c}(g)\right)^{\star}
\;\; \Rightarrow \;\; \alpha
\in
\left((a_{1}, a_{2}, \cdots , a_{n})D\right)^\star  =J^\star \,.
\]

Therefore $J^{\star_{_{\!{\ell}}}} \subseteq J^\star$. 

 Finally, since $\Na(D, \star) \subseteq \KN(D, \star)$ (Theorem 
\ref{th:1.21} (7)) then, for each $E \in
\overline{\boldsymbol{F}}(D)$, 
$E^{\tilde{\star}} = E\Na(D, \star)  \cap K \subseteq E \KN(D, 
\star)  \cap K = E^{\star_{_{\!\ell}}}$.   
\hfill $\Box$

\bigskip

 The statements (1), (2) and (3) of our next result are 
essentially  the ``$\KN$--analogues'' of Proposition \ref{pr:1.12} 
and Corollary \ref{co:1.16};  the statements (4) is an 
``$\calL$--analogue'' of Lemma \ref{lm:1.1} (9).

%COROLLARY 6.4
\begin{coor} \label{co:1.23} \sl
Let $\star$ be a semistar operation on an
             integral domain $D$ and let ${\calL} ={\calL}(D, \star)$
be the set of all $\star$--monolocalities of $D$. Then:
               \bara

% 1
\bf \item \sl $
             \Na(D, \star) = \Na(D, \widetilde{\star}) = \KN(D,
             \widetilde{\star})\,,$  (in particular, 
$\widetilde{\star} =(\widetilde{\star})_\ell$).  
%2
\bf \item \sl  $ \KN(D,  \star_{a}) =\Kr(D, \star_{a}) = \Kr(D,
\star)\,,$  (in particular, $\star_{a} =(\star_{a})_\ell$).  

%3
\bf \item \sl  $D$ is a P$\star$MD if and only if
$ \KN(D,  \widetilde{\star})=\KN(D,  \star_{a})$.
%4
    \bf \item \sl   For each $(L,N) \in \calL$, $N\cap D$ is a
quasi--$\widetilde{\star}$--prime  of $D$.
    \eara
\end{coor}

    \bf Proof. \rm  \bf (1) \rm is an easy consequence of Theorem
\ref{th:1.21} (7) and Proposition \ref{pr:1.12} (1). 

    \bf (2) \rm follows from the fact that, for the e.a.b. semistar
operation $\star_a$, each
    $F\in \boldsymbol{f} (D)$ is $\star_a$--e.a.b., hence $\calL(D,
\star_a)$ is the set of all the $\star$--valuation
    overrings of $D$   (since $T \in \calL(D,
\star_a)$ is necessarily a valuation overring of $D$ and $T = 
T^{\star_a}$, i.e. $T$ is a $\star_a$--valuation overring (Lemma 
\ref{lm:1.1} (6)),  or equivalently  a $\star$--valuation overring, 
of $D$ \cite[Proposition 3.3]{FL2}).  

    \bf (3) \rm  follows from (1), (2) and  \cite[Theorem 3.1 and 
Remark
3.1]{FJS}.  .

    \bf (4) \rm Using Theorem \ref{th:1.21} (6) and (7), we have
$(N\cap D)^{\widetilde{\star}} =
    (N\cap D)\Na(D, \star) \cap K \subseteq (N\cap D)\KN(D, \star)
\cap K  \subseteq (N(X) \cap  \KN(D, \star)) \cap K$.
    Therefore $(N\cap D)^{\widetilde{\star}} \cap D \subseteq ((N(X)
\cap  \KN(D, \star)) \cap K) \cap D =
    ((N(X) \cap  \KN(D, \star)) \cap L)  \cap D = N \cap D$\,.
\hfill  $ \Box$

    %REMARK 6.5
  \begin{reem} \label{rk:6.5}  \rm So far we have given no indication 
that
  $\KN(D, \star)$ and $\SKN(D, \star)$ are ever different.
  In Example \ref{ex:1.26} we exhibit a
   (semi)star operation (of finite type) $\star$ on a Noetherian
   integrally closed integral domain $D$ such that
   $\KN(D, \star) \subsetneq  \SKN(D, \star)$  and thus, in 
particular,
   ${\calL}(D, \star)  \supsetneq    {\calL'}(D, \star)$  (and so 
$\wedge_{\calL}=\star_{_{\!{\ell}}} \lneq {\star_{_{\!{\ell'}}}}\leq 
\wedge_{{\calL}'}$).   
   Moreover in this example we will see that:
   $$      \begin{array}{rl}
   (\star =) \ \star_{_{\!f}}   \lneq  & \hskip -8pt
   {\star_{_{\!{\ell'}}}} \ (=\star_{a} =t_{D}) \,,\\
  (d_{D}=) \ \widetilde\star \
   =& \hskip -8pt {\star_{_{\!{\ell}}}}   \lneq \star_{_{\!f}} \,.
 \end{array}
 $$
       
     Moreover, it is not difficult to give an example of a
         semistar operation $\star$ such that
           ${\widetilde{\star}}   \lneq
           {\star_{_{\!{\ell}}}}$ (cf. the following Example 
\ref{ex:7.5}).

       \end{reem}
      \vskip 12pt

Now that we have made note that $\calL$ and $\calL'$
are not always the same we investigate the implications of assuming
that they or the related rational function rings or the related
semistar operations are the same.

    %
    %%COROLLARY 6.6
 
              \begin{coor} \label{cor:1.22}  \sl  Let $\star$ be a
semistar operation on an
integral domain $D$.  Then the following are equivalent.

\brom
\bf \item[(i)]  \sl ${\calL} =  {\calL'}$\,;
\bf \item[(ii)]  \sl   $\wedge_{\calL} =\wedge_{\calL'}$\,;
\bf \item[(iii)] \sl $\star_{_{\!f}}  \leq  {\star_{_{\!{\ell}}}} $\,;
\bf \item[(iv)] \sl ${\star_{_{\!{\ell}}}} = \star_{_{\!f}}$\,;
\bf \item[(v)] \sl  $ {\star_{_{\!{\ell}}}} =\wedge_{\calL}
=\wedge_{\calL'} =
{\star_{_{\ell'}}} = \star_{_{\!f}}$\,.
\erom
\end{coor}

\bf Proof.  \rm  (v) $\Rightarrow$   (iv) $\Rightarrow$ 
(iii) and (i)
$\Rightarrow$
(ii) are trivial.

(iii) $\Rightarrow$ (i) because  (iii) implies that each $L \in
             {\calL}(D, \star)$ is a $\star$--overring of
             $D$, since ${\star_{_{\!{\ell}}}} = \wedge_{\calL}$.

(ii) $\Rightarrow$ (iv) because we know that, in general,
         $\wedge_{\calL} \leq \star_{_{\!f}} \leq \wedge_{\calL'}$
 (Propositions \ref{pr:6.2} and \ref{pro:6.3}).   

 (i) $\Rightarrow$ (v) is obvious, using the fact that we already 
know that (i) $\Leftrightarrow$ (iv).   
\hfill $\Box$

% COROLLARY 6.7
\begin{coor} \label{cor:6.7}  \sl  Let $\star$ be a
semistar operation on an
integral domain $D$.  If  $\star_{_{\!f}}$ is stable (i.e.
$\widetilde{\star}= \star_{_{\!f}}$  \cite[Corollary 3.9]{FH}) 
   then $\KN(D, \star) = \SKN(D,
\star)$ and $ {\star_{_{\!{\ell}}}} = \wedge_{\calL} =\wedge_{\calL'}
=
{\star_{_{\ell'}}} \ (= \star_{_{\!f}}=\widetilde{\star}).$
 \end{coor}

\bf Proof. 
\rm  The result follows from   Corollary \ref{cor:1.22} and   the 
fact that, in general, $ 
\widetilde{\star} \leq {\star_{_{\ell}}} \leq \star_{_{\!f}}$ 
(Proposition \ref{pro:6.3}).
  \hfill $\Box$
  
             \hskip 12pt

%%%%%%%%%%%%%%%%%%%%%    SECTION 7

\section{Constructions and Examples}

Thus far, we have focussed on the situation where we begin with a
semistar operation and investigate related overrings determined by
the e.a.b.--ideals associated to the semistar operation.  Now we
reverse that and begin with a collection of overrings and use them
to defined a semistar operation.  The major questions will concern
how the (strong) monolocalities relate to the defining
collection of overrings.

% PROPOSITION 7.1
\begin{prro} \label{pr:1.25}  \sl Let $\calT :=
\{T_\lambda \mid \lambda \in \Lambda\}$ be a family of  quasilocal 
 overrings
of an integral domain $D$ and let $F \in \boldsymbol{f}(D)$. Then
$F$ is $\wedge_{\calT}$--e.a.b. if and only if $FT_{\lambda}$ is
principal as a fractional ideal of $T_{\lambda}$, for each $\lambda
\in
\Lambda$\,.
\end{prro}

\bf Proof.  \rm The ``only if'' part.

\bf Claim. \sl Let  $\lambda \in
\Lambda$ and  let $F, G, H \in  
\boldsymbol{f}(D)$. Assume that $F$ is
$\wedge_{\calT}$--e.a.b. and   $FGT_{\lambda}\subseteq
FHT_{\lambda}$\,, \
then $GT_{\lambda}\subseteq HT_{\lambda}$.\rm

Note that $FGT_{\lambda}\subseteq FHT_{\lambda}$ implies that
$(FG)^{\wedge_{\calT}} \subseteq
(FGT_{\lambda})^{\wedge_{\calT}}\subseteq
(FHT_{\lambda})^{\wedge_{\calT}}$.
 Since  $F$ is
$\wedge_{\calT}$--e.a.b., thus also
$(\wedge_{\calT})_{_{\!f}}$--a.b. (Lemma \ref{le:1.5}  (2)),    
then
 $G^{(\wedge_{\calT})_{_{\!f}}} = G^{\wedge_{\calT}}
\subseteq ({HT_{\lambda}})^{(\wedge_{\calT})_{_{\!f}}} \subseteq
({HT_{\lambda}})^{\wedge_{\calT}}
$\,. \ Therefore  $ GT_{\lambda}=(G^{\wedge_{\calT}})T_{\lambda}
\subseteq
(({HT_{\lambda}})^{\wedge_{\calT}})T_{\lambda}  =HT_{\lambda}$\
(Lemma
\ref{lm:1.1} (10))\,.

The conclusion of the ``only if'' part follows from   the previous 
Claim and   Remark  
\ref{rk:3.5},   since $T_{\lambda}$
is  quasilocal   and each finitely
generated $T_{\lambda}$--submodule of $K$, $G_{\lambda} \in
\boldsymbol{f}(T_{\lambda})$,  is of the type
$GT_{\lambda}$, for some $G \in \boldsymbol{f}(D)$.

For the ``if'' part, assume that $F, G, H \in \boldsymbol{f}(D)$,
$(FG)^{\wedge_{\calT}} \subseteq
(FH)^{\wedge_{\calT}}$ and
$FT_{\lambda}$ is
principal as a fractional ideal of $T_{\lambda}$, for each $\lambda
\in
\Lambda$\,. Then, clearly, $FGT_{\lambda}=
(FG)^{\wedge_{\calT}}T_{\lambda} \subseteq
(FH)^{\wedge_{\calT}}T_{\lambda} =FHT_{\lambda}$. Since
$FT_{\lambda}$ is principal, then $GT_{\lambda}
\subseteq HT_{\lambda}$, for each $\lambda \in
\Lambda$\,, \ and thus $G^{\wedge_{\calT}} \subseteq
H^{\wedge_{\calT}}$\,. \hfill $\Box$

\bigskip

We digress momentarily to give a corollary to the last result which
illustrates some nice closure properties of the $(-)_{\ell}$ 
operation.

%COROLLARY 7.2
\begin{coor} \label{cor:1.26} \sl
 Let $\star$ be a semistar operation on an
integral domain $D$.
\bara
\bf \item \sl Let $F \in \boldsymbol{f}(D)$. If $F$ is
$\star$--e.a.b. then $F$ is
$\star_{\ell}$--e.a.b..
\bf \item \sl $\calL(D, \star) \subseteq \calL(D, \star_{\ell})$\
(more precisely, $\calL(D, \star) = \{ L \in \calL(D, \star_{\ell})
\mid L = L^{\star_{_{\!f}}} \}$).

\bf \item \sl $\KN(D, \star) = \KN(D, \star_{\ell})$\,.

\bf \item \sl $(\star_{\ell})_{{\ell}}= \star_{\ell}$\,.
\eara
\end{coor}
\bf Proof. \rm Recall that ${\star_{\ell}} = {\wedge_{\calL}}$,
where $\calL = \calL(D, \star)$  (Proposition \ref{pro:6.3}). 

\bf (1) \rm Assume that $F$ is
$\star$--e.a.b. then, by definition of $\star$--monolocality, $FL$ is
principal for each $L \in
\calL(D, \star)$. Therefore, by Proposition \ref{pr:1.25}, $F$ is
${\wedge_{\calL}}$--e.a.b. \ (= $\star_{\ell}$--e.a.b.).

\bf (2)  \rm Let $L \in \calL(D, \star)$ and let $F$ be
${\star_{\ell}}$--e.a.b. \ (= ${\wedge_{\calL}}$--e.a.b.). As above,
by Proposition \ref{pr:1.25},
we know that $FL$ is principal, thus $L$ belongs also to  $\calL(D,
{\star_{\ell}})$, since $L = L^{\star_{_{\!f}}}$ and $\star_{\ell}
\leq {\star_{_{\!f}}}$  (Proposition \ref{pro:6.3}). 

 For the parenthetical statement,   let $L \in \calL(D, 
\star_{\ell})$ and let $F$ be
$\star$--e.a.b.. By (1)  $F$ is
${\star_{\ell}}$--e.a.b. \ (= ${\wedge_{\calL}}$--e.a.b.), thus $FL$
is principal,
for each $L \in \calL(D, \star)$ (Proposition \ref{pr:1.25}).
We conclude that $L \in
\calL(D, \star_{\ell})$ belongs to $\calL(D, \star)$ if and only if
$ L = L^{\star_{_{\!f}}}$.

\bf (3) \rm From (2) we deduce that $\KN(D, \star) \supseteq \KN(D,
\star_{\ell})$. The opposite inclusion follows from Theorem
\ref{th:1.21} (1) and (7), since if $g \in D[X]$ is such that
$\boldsymbol{c}(g)$ is $\star$--e.a.b.
then $\boldsymbol{c}(g)$ is also $\star_{\ell}$--e.a.b. by (1).

\bf (4)  \rm is a straightforward consequence of (3). \hfill $\Box$

\medskip 
We now turn back to the special case where we begin with a collection
of overrings of a domain $D$ and use them to define a semistar
operation.

         %
         %EXAMPLE 7.3
\begin{prro} \label{ex:1.27}
    \sl  Let $\calT := \{T_\lambda \mid \lambda \in \Lambda\}$ be a
family of   quasilocal   overrings of an integral domain $D$
    and set $\ast := \wedge_{\calT}$. Then:
    \begin{enumerate}
    %a
      \rm \bf \item[(1)] \sl   $\calT \subseteq \calL' (D, \ast) 
\subseteq
\calL (D, \ast)$\,.

        %b
      \rm \bf \item[(2)] \sl   $
    \wedge_{\calL} \leq \wedge_{\calL'} \leq  \ast
    $ \ and \ $  \left(\wedge_{\calL'}\right)_{_{\!f}}
=\ast_{_{\!f}}$\,.

    %c
       \rm \bf \item[(3)]  \sl   $ \KN(D, \wedge_{\calL}) =\KN(D,
\wedge_{\calL'}) =
    \KN(D, \ast)$\;  and \; $ \KN(D, \ast) =\SKN (D, \ast) $\,.

    %d
      \rm \bf \item[(4)]  \sl     ${\ast_{_{\!{\ell}}}} =
{\ast_{_{\!{\ell'}}}} =
    \wedge_{\calL} \leq (\wedge_{\calL'})_{_{\!f}} =\ast_{_{\!f}}$\,
    and \,  ${\ast_{_{\!{\ell}}}} =
({\ast_{_{\!{\ell'}}}})_{_{\!{\ell}}} =
    (\wedge_{\calL})_{_{\!{\ell}}} =
    (\wedge_{\calL'})_{_{\!{\ell}}}$\,.
    \end{enumerate} \rm 
\end{prro}

    \medskip 

    \bf Proof.   \bf  (1) \rm    Each $T_{\lambda}$ is obviously a 
$\ast$--overring
of $D$,
    since if $F \in \boldsymbol{f}(D)$ then  $FT_{\lambda} = F^\ast
T_{\lambda} =
    (FT_{\lambda})^\ast$ (Lemma \ref{lm:1.1} (10)).
    Furthermore, note that if $F \in \boldsymbol{f}(D)$ is
$\ast$--e.a.b.
    then $FT_{\lambda}$ is principal in  $T_{\lambda}$, for each
$\lambda \in \Lambda$.
    As a matter of fact, if $G, H \in \boldsymbol{f}(D)$ are such
that $FGT_{\lambda} \subseteq  FHT_{\lambda}$,
    then $ (FGT_{\lambda})^\ast =  FGT_{\lambda} \subseteq
FHT_{\lambda} =(FHT_{\lambda})^\ast$, thus
    $ (GT_{\lambda})^\ast =  GT_{\lambda} \subseteq  HT_{\lambda}
=(HT_{\lambda})^\ast$, because $F$ is $\ast$--e.a.b..
    Therefore $FT_{\lambda}$ is quasi-cancellative and so it is
principal in $T_\lambda$  (Remark \ref{rk:3.5}).  

      \bf (2) \rm     From   (1),   we deduce immediately that:
    $$
    \wedge_{\calL} \leq \wedge_{\calL'} \leq \wedge_{\calT} = \ast\,.
    $$
    Moreover, $\ast_{_{\!f}} =
\left(\wedge_{\calL'}\right)_{_{\!f}}$, since
    for each $F \in \boldsymbol{f}(D)$, we have:
    $$
    \begin{array}{rl}
    F^{\wedge_{\calL'} } =& \hskip -5pt ( \bigcap \{FT_{\lambda} \mid
\lambda \in \Lambda \}) \cap
    ( \bigcap \{ FL' \mid L'  \in { \calL'}\setminus {\calT}  \}) \\
    = & \hskip -5pt
    F^\ast  \cap  ( \bigcap \{ FL' \mid L'  \in { \calL'}\setminus
{\calT}
    \}) \\
   = & \hskip -5pt
    F^\ast  \cap  ( \bigcap \{  F^{\ast}L'   \mid L'  \in {
\calL'}\setminus {\calT}
    \})  = F^\ast \,.
    \end{array}
    $$

      \bf (3) \rm     Since $ \KN(D, \ast) =\KN(D, \ast_{\ell}) $
    (Corollary \ref{cor:1.26} (3)),  ${\ast_{\ell}}   =
{\wedge_{\calL}}$
    (Proposition \ref{pro:6.3}),   and $\ast_{_{\!f}} =
    \left(\wedge_{\calL'}\right)_{_{\!f}}$ by (b),
    then we easily deduce  that
    $ \KN(D, \wedge_{\calL})=  \KN(D, \ast_{\ell})    = \KN(D, 
\ast)
=\KN(D, \wedge_{\calL'}) $\,.

    From Lemma
\ref{lm:1.1} (10) we deduce immediately:

    \bf Claim 1. \sl Let $g \in D[X]$,\ $g \neq 0$, then
    $\boldsymbol{c}(g)^{\wedge_{\calT}} = D^{\wedge_{\calT}}$ if and
    only if $\boldsymbol{c}(g)T_{\lambda} = T_{\lambda}$, for each
    $\lambda \in \Lambda$\,. \rm  

\smallskip

  Claim 1 implies: 

    \bf Claim 2. \sl $\Na(D, {\wedge_{\calT}}) = \bigcap \{\Na(D,
    \star_{\{T_{\lambda} \}}) \mid \lambda \in \Lambda\} =\bigcap
    \{T_{\lambda}(X) \mid \lambda \in \Lambda\}\,.$ \rm
    
    \smallskip

    From Claim 2 and from the fact that $\calT \subseteq \calL'(D,
    \ast)$  we deduce that $\SKN(D, \ast) =\bigcap
    \{L'(X) \mid L' \in \calL'(D, \ast)\} \subseteq \bigcap
    \{T_{\lambda}(X) \mid \lambda \in \Lambda\} =\Na(D,
    \wedge_{\calT}) \subseteq  \KN(D,
    \wedge_{\calT}) =\KN(D,
    \ast) $.  Since in general, $\KN(D,
    \ast)  \subseteq \SKN(D, \ast)$, we conclude that $\KN(D,
    \ast)  = \SKN(D, \ast)$\,.

      \bf (4) \rm      From   (3)   and from  Proposition 
\ref{pro:6.3}, 
  we
have that
    $ {\ast_{_{\!{\ell}}}}    = \wedge_{\calL} =  
{\ast_{_{\!{\ell'}}}}$   
    and so, by Corollary \ref{cor:1.26}   (4),  
${\ast_{_{\!{\ell}}}}   =
(\wedge_{\calL})_{_{\!{\ell}}}
=( {\ast_{_{\!{\ell'}}}})_{_{\!{\ell}}}$   
    From (b)  we know that $\left(\wedge_{\calL'}\right)_{_{\!f}}  =
    \ast_{_{\!f}}$. Since
    $\wedge_{\calL} \ (=  {\ast_{_{\!{\ell}}}})$    is a semistar
    operation of finite type then, clearly,  $\wedge_{\calL} \leq
    \left(\wedge_{\calL'}\right)_{_{\!f}}$. Finally, note that
$(\left(\wedge_{\calL'}\right)_{_{\!f}})_{_{\!{\ell}}} =
\left(\wedge_{\calL'}\right)_{_{\!{\ell}}}$ and
$( \ast_{_{\!f}})_{_{\!{\ell}}}   = \ast_{_{\!{\ell}}}$.
\hfill $\Box$

\vspace{.2in}

If we assume in addition to the hypotheses of Proposition 
\ref{ex:1.27} 
that each $T_{\lambda}$ is integrally closed we can prove a little 
more.

% COROLLARY 7.4
\begin{coor} \label{co:1.28} Suppose in addition to the hypotheses of 
Proposition \ref{ex:1.27} that each $T_{\lambda}$ is integrally 
closed.  
Then  $\KN(D, \ast) =\SKN (D, \ast)  =    \bigcap  \{T_{\lambda}(X) 
\mid \lambda \in \Lambda\}\,.$  
\end{coor}

\medskip

{\bf Proof.}  We already know that 
$\KN(D,\ast) = \SKN(D,\ast)$.  Since each $T_{\lambda}$ is a 
strong $\ast$-monolocality    (Proposition \ref{ex:1.27} (1)),   it 
follows immediately from the 
definitions that $\SKN(D,\star) \subseteq \bigcap \{T_{\lambda}(X) 
\mid \lambda \in \Lambda\}$.  We 
will finish the proof by showing that 
$\bigcap\{T_{\lambda}(X) \mid \lambda \in \Lambda\} \subseteq 
\KN(D,\ast)$.  Let 
$f/g\in \bigcap\{T_{\lambda}(X) \mid \lambda \in \Lambda\} $ and 
suppose that 
$f,g \in D[X]$,   with $g \neq 0$,   and that 
they have no common factors over $K[X]$.  Choose a particular 
$T_{\lambda}$. We will consider content ideals of polynomials as 
ideals of $T_{\lambda}$.  By definition, we 
know that $f/g = h/k$ where $h,k \in T_{\lambda}[X]$ and 
$\boldsymbol{c}_{T_{\lambda}}(k) = T_{\lambda}$.   Since $f, g$ have 
no common factors over $K[X]$,    by Euclid's Lemma    
we know that $g$ is a factor of $k$    in $K[X]$.    If we 
rewrite $h/k$ as $dh/dk$ for an appropriate   nonzero   constant 
$d \in T_{\lambda}$,  we have $g$ being a factor over 
$T_{\lambda}[X]$ 
of $dk$ and $\boldsymbol{c}_{T_{\lambda}}(dk) = dT_{\lambda}$.  Then 
since $T_{\lambda}$ is 
integrally closed we know that $g$ has invertible (hence 
principal) content    \cite[Theorem 1.5]{MNZ}.     Finally, since 
$\boldsymbol{c}_{T_{\lambda}}(h) 
\subseteq\boldsymbol{c}_{T_{\lambda}}(k)= T_{\lambda}$ and 
$f/g = h/k$ it follows easily that $\boldsymbol{c}_{T_{\lambda}}(f) 
\subseteq \boldsymbol{c}_{T_{\lambda}}(g)$. 

 Since we 
have been working with an arbitrary $T_{\lambda}$ we have proven:

\begin{itemize}
\item $\boldsymbol{c}_D(f)^{\ast} \subseteq 
\boldsymbol{c}_D(g)^{\ast}$ (with content ideals 
considered now as ideals of $D$).

\item $\boldsymbol{c}_D(g)$ is a $\ast$-e.a.b. ideal of $D$ since its 
extension to 
each $T_{\lambda}$ is principal.

\end{itemize}

This proves that $f/g \in \KN(D,\ast)$ which finishes the proof. 
\hfill $\Box$

\medskip 

\rm  In the setting where we begin with a
collection of overrings of a domain $D$ it can seem that differences
between strong monolocalities and monolocalities and the associated
constructions should disappear.  In particular, we could hope that
the inequality in part    {(4)}   of   Proposition \ref{ex:1.27}   
should be an equality.
We next give an example to demonstrate that this inequality can be
strict.

%% EXAMPLE 7.5

\begin{exxe} \label{ex:1.28}
\rm Let $k$ be a field and    let 
$\{X,Y,W,X'_{1},Y'_{1},X'_{2},Y'_{2}, \ldots; Z\}$ be an infinite 
family
of indeterminates over $k$. Set    $R := 
k[X,Y,W,X'_{1},Y'_{1},X'_{2},Y'_{2}, \ldots ]$.    Let $M$ be 
the maximal ideal of $R$ generated by    
$\{X,Y,W,X'_{1},Y'_{1},X'_{2},Y'_{2}, \ldots \}$,    let $D := R_{M}$ 
and let $K$ be the quotient field of $D$.  
  For each positive integer $i$ define 
$$
T_{i} : =    D\left[\frac{W}{XX'_{i}+YY'_{i}}\right]\,.   
$$
Let $\calT := \{T_i \mid i > 0\}$ and set  $\ast:= \wedge_{\calT}$ 
(i.e.  $E^{\ast} := \bigcap_{i \geq 1}  ET_{i}$, for each $E\in 
\overline{\boldsymbol{F}}(D)$).
Also,  let $P_{W}$ be the prime ideal of $D$ generated by 
$W$.  Note that $\overline{k} := D_{P_{W}}/P_WD_{P_{W}}$ is isomorphic to 
$k(X,Y,X'_{1},Y'_{1},X'_{2},Y'_{2}, \ldots )$.  Let $\varphi$ be the 
canonical homomorphism from $D_{P_{W}}$ to $\overline{k} = 
D_{P_{W}}/P_WD_{P_{W}}$.  Let $\overline{V}$ be a minimal valuation 
overring of $\overline{D}:= D/P_W$ (in its quotient field isomorphic to $\overline{k}$) and let $V = \varphi^{-1}(\overline{V})$.  
Then $V$ is a  
minimal valuation overring of $D$ which has 
$D_{P_{W}}$ as an overring.

\begin{enumerate}

% 1
\bf \item  \sl $V$ is a minimal valuation overring of each $T_{i}$. 
\rm

%2
\bf \item  \sl $V(Z)$ is a minimal valuation overring of each 
$T_{i}(Z)$. \rm

%3
\bf \item  \sl  $\KN(D,\ast) = \SKN(D,\ast) = \bigcap_{i \geq 1}  
T_{i}(Z)$

%4
\bf \item  \sl  Let $M_{V}$ be the contraction of the maximal ideal 
of \
$V(Z)$ to $\KN(D,\ast)$.  Then 
$M_{V}$ is a maximal ideal of
$\KN(D,\ast)$.

%5
\bf \item  \sl $M_{V}$ is the only maximal ideal of 
$\KN(D,\ast) $.

%6
\bf \item  \sl   $D^{\ast} = D$.

 %7
\bf \item  \sl $(\KN(D,\ast) = \SKN(D,\ast) =) \  \bigcap_{i \geq 1} 
T_{i}(Z) = D(Z)$.

%8
\bf \item  \sl  $D$ is a $\ast$--monolocality but not a strong 
$\ast$--monolocality.

%9 
\bf \item  \sl  The (semi)star operation $\ast$ is such that the 
inequality in Proposition \ref{ex:1.27} (4) is 
strict (i.e., 
$\wedge_{\calL(\ast)} \lneq (\wedge_{\calL'(\ast)})_{_{\!f}} $).   
\end{enumerate}

 %%%%%%%%%%%%%%%%
\smallskip

  {\bf Proof.} \bf  (1) \rm  is a consequence of the fact that if 
$\frac{W}{XX'_{i}+YY'_{i}} \notin V$ then $\frac{XX'_{i}+YY'_{i}}{W} 
\in V \subseteq D_{P_W}$, which is a contradiction.

\bf (2) \rm  follows from (1) and \bf (3) \rm  is a consequence of
 Corollary \ref{co:1.28}, since each $T_i$ is integrally closed.  
(The claim that each $T_i$ is integrally closed follows easily from 
\cite[Theorem 2]{Sally}.)

\bf (4) \rm  Proposition \ref{ex:1.27} implies that each $T_{i}$ is a 
$\ast$-monolocality.  It follows then from Theorem \ref{th:1.21} that 
$T_{i}(Z)$ is an overring of $\KN(D,\ast)$ for each $i$.  Hence 
$V(Z)$ is an overring of $\KN(D,\ast)$
We assumed $V$ to be a minimal valuation overring of $D$.  It follows 
that $V(Z)$ is a minimal valuation overing of $D(Z)$.  It is clear 
then that $V(Z)$ is also a minimal valuation overring of any ring 
properly in between $D(Z)$ and $V(Z)$.  In particular, $V(Z)$ is a 
minimal valuation overring of $\KN(D,\ast)$.  The result follows 
immediately.

\bf (5) \rm  Let $d \in M_{V}$.  Then $d \in \KN(D,\ast)$.  As we 
noted above, $T_{i}(Z)$ is an overring of $\KN(D,\ast)$.  Also recall 
that    $V(Z)$    is a minimal valuation overring of each    
$T_i(Z)$.     In particular, $d$ is a nonunit in each   $T_i(Z)$.    
Hence,   
$1+d$ is a unit in each 
$T_{i}(Z)$   and so, by (3), it is a unit in $\KN(D,\ast)$.

\bf (6) \rm  It is easy to see that $\bigcap_{i\geq 1} T_{i} = D$,   
hence $D^\ast = \bigcap_{i\geq 1} DT_{i}= \bigcap_{i\geq 1} T_{i} =D$ 
.

\bf (7) \rm  Combine Theorem \ref{th:1.21} (6) and (7) with the fact 
that we already know that $\KN(D,\ast)$ is quasilocal   and that 
$\KN(D,\ast) \cap K= \bigcap_{i\geq 1} T_{i}(Z) \cap K 
=\bigcap_{i\geq 1} T_{i} = D$.  

\bf (8) \rm $D$ being a $\ast$--monolocality follows  from the 
fact that $\KN(D,\ast) = D(Z)$ (Theorem \ref{th:1.21}  (5)).  $D$ is 
not a strong $\ast$--monolocality 
  because it is not a $\ast$--overring: let $J = (X,Y)D$, then 
$W \in J^{\ast}$ but $W \not \in JD=J$.

  \bf (9)  \rm  Since $D$ itself is a 
$\ast$--monolocality (by (8)), then 
$\wedge_{\calL}$ is the identity function.  However, 
the proof of (8) above demonstrates that  
$\ast_f \ ( = (\wedge_{\calL'})_{_{\!f}}$ \  by Proposition 
\ref{ex:1.27} (4)) is not the identity.

\end{exxe}

\medskip

So we have an example of a semistar operation $\star$ on a domain $D$ 
of the form  $\wedge_{\calT}$, derived from    a family $\calT$ of   
overrings of $D$, such that 
$\KN(D,\star) = \SKN(D,\star)$ but there exists a 
$\star$--monolocality which is not strong   (Example 
\ref{ex:1.28}).

We next give an example of a semistar operation,   also defined  by 
a  $\wedge$--construction,    in which 
things work exactly as one might hope.

% EXAMPLE 7.6

\begin{exxe} \label{ex:11.43} \rm 
Let $k$ be a field and let $R :=k[X,Y]$ be the ring of polynomials 
over $k$ in the two variables $X$ and $Y$.  Let $P := (X,Y)$ be the 
maximal ideal of $R$ generated by the variables and let 
$D := R_{P}$.  Let $M$ denote the maximal ideal of the local ring 
$D$.  Set

\begin{itemize}

\item $D_{1} := D[X/Y]$\,,\;\;\; $D_{2} := D[Y/X]$\,;

\item  $\calbT:=\{T_\lambda \mid T_\lambda\in \Lambda\}$ is the 
collection of all localizations of $D_{1}$ and $D_{2}$ at their 
maximal ideals.

\end{itemize}

  \sl Set
$\ast:= \wedge_{\calT}$. Then $\ast$ is an example of a semistar 
operation on a  integrally closed Noetherian local domain $D$  such 
that:
 
\bara 

\bf \item  $ {\ast_{_{\!{\ell}}}} =\wedge_{\calL} 
=\wedge_{\calL'} = {\ast_{_{\ell'}}} = \ast_{_{\!f}} = \ast$\;  \sl 
(i.e. $\calL = \calL'$, by Corollary \ref{cor:1.22}).
\bf \item  $ \widetilde{\ast} \lneq \ast \lneq \ast_{a}\,.$
\eara \rm    

We know by Proposition \ref{ex:1.27} (1) 
that $\calT \subseteq \calL' (D,   \ast)    \subseteq
\calL (D,   \ast)$.    Suppose then that $T$ is a 
$\ast$--monolocality 
of $D$.  It is clear from Proposition \ref{pr:1.25} that the ideal 
$I = (X,Y)$ is a   $\ast$--e.a.b. ideal of $D$.    Hence the ideal 
$I$ must extend to a 
principal ideal in $T$.  Since $T$ is quasilocal it follows that $IT$ 
is generated by either $X$ or $Y$.  Hence either $Y/X$ or $X/Y$ lies 
in $T$.  Hence $T_{\lambda} \subseteq T$ for some 
$T_{\lambda} \in \calT$.  It follows that $\calT$ consists exactly of 
the minimal   $\ast$--monolocalities    of $D$.  It follows from this 
then 
that the   $\ast$--monolocalities    and the strong   
$\ast$--monolocalities  
coincide.  This is sufficient to prove that 
$$
{\ast_{_{\!{\ell}}}} =\wedge_{\calL} 
=\wedge_{\calL'} = {\ast_{_{\ell'}}} = \ast_{_{\!f}} = \ast\,.$$

Now observe that the 
$T_{\lambda}$'s are neither localizations of $D$ nor valuation 
overrings of $D$.  This proves that $ \widetilde{\ast} \lneq \ast 
\lneq \ast_{a}$.  Suppose that, in fact, 
$\widetilde{\ast} = \ast$.  Then Corollary \ref{co:1.23} indicates 
that 
$\KN(D,\ast) = \Na(D,\ast)$.  Recall that the localizations of 
$\Na(D,\ast)$ at maximal ideals have the form 
$D_P(X)$ where $P$ is a prime ideal of $D$.  Similarly, if $\ast = 
\ast_a$ then Corollary \ref{co:1.23} indicates that 
$\KN(D,\ast) = \Kr(D,\ast)$.  Recall that the localizations of   
$\Kr(D,\ast)$    at maximal ideals have the form $V(X)$ where $V$ is 
a valuation overring of $D$.  In the present setting the 
localizations of $\KN(D,\star)$ are exactly the rings $T_i(X)$.  The 
result follows immediately.

\end{exxe}

Example \ref{ex:11.43} is significant because we indicated that an 
important objective of this 
article was to demonstrate that the Nagata ring construction and the 
Kronecker function ring construction were at opposite ends of a 
spectrum.  For the generalization to have any real power we need to 
demonstrate that we can find something which is properly in between 
these two extremes.  We also indicated that we wanted to give a 
method 
for approximating a given semistar operation by a semistar operation 
which was constructed by means of extension to a collection of 
overrings.  We have given such a mechanism, but again we need to show 
that this is meaningful by demonstrating that the semistar operation 
obtained can turn out to be associated with a collection of overrings 
which consists neither of localizations nor of valuation overrings of 
the domain $D$.  In this example we have 
 a star operation    $\ast$    such that 
   $\ast$   is equal to all four of the approximations developed in 
this work    (Definition \ref{def:6.1} ).    And yet we also have    
\[
 \widetilde{\ast } \lneq \ast\lneq \ast _{a}\,.
  \]   
This indicates that    $\ast$   and all of its KN and $\SKN$ 
derivatives lie 
properly in between   `` the localization constructions'' associated 
to 
  $\widetilde{\ast}$   and  ``the valuation domain constructions'' 
associated to   $\ast_{a}$.

    \medskip

%%%%%%

  %EXAMPLE 7.7
\begin{exxe} \label{ex:1.26}
 \sl Example of a (semi)star operation $\star$ on a Noetherian
 integrally closed    domain $D$ such that  $ 
\star_{_{\!{\ell}}}  
             \lneq \star \ (=  \star_{_{\!f}})  \lneq
             {\star_{_{\!{\ell'}}}}$ and so   $\SKN(D, \star) 
\subsetneq \KN(D,
\star)$,   in particular    $\calL(D, \star) \subsetneq 
\calL'(D, \star)$. 

\smallskip
\rm
Let $k$ be a field, $D:= k[X,Y]_{(X,Y)}$   $M :=(X,Y)D$ and let
$K:=k(X,Y)$\,.
Using the techniques of \cite[(32.4)]{G}, we construct a new
(semi)star
  operation $\star$ on $D$ as follows:

  \begin{enumerate}
%1
  \rm \bf \item[1.] \rm If $\,dD\,$ is any nonzero principal ideal of
  $\,D\,,$ \ then
  $\,(dD)^\star := dD\,$.

%2
  \rm \bf \item[2.]   \rm  If $\,J \subseteq D\,$ is a nonzero ideal
of
  $\,D\,$
  which is not contained in any proper principal ideal of $\,D\,$, \
  then
  $\,J^\star :=   M\,$.   

%3
  \rm \bf \item[3.]  \rm If $\,J \subseteq D\,$ is a nonzero ideal of
  $\,D\,$
  which is not principal,
   but is contained in a principal ideal, then we factor $\,J\,$ as
  $\,J =
  \alpha I\,,$ \ where
  $\,\alpha\,$ is a GCD of a set of generators of $\,J\,$ and $\,I :=
(J
  :_{D} \alpha D)\,$ is not contained
  in any proper principal ideal of $\,D\,$\ by the choice of
$\,\alpha \,.$ \
  Then
  $\,J^\star := \alpha  M\,$.

%4
  \rm \bf \item[4.]   \rm  If $\,J\,$ is a nonzero fractional ideal
of
  $\,D\,$
  which is
  not contained in $\,D\,$,\ choose a nonzero element $\,\beta\in
D\,$ such
  that $\,\beta J \subseteq D\,$. \
  Then define
  $\,J^\star := (1/\beta )(\beta J)^\star\,$.

% 5
  \rm \bf \item[5.]  \rm If $\,J \in \overline{\boldsymbol{F}}(D)
\setminus
  \boldsymbol{F}(D)\,$ we define $\,J^\star :=  K \,$.   

  \end{enumerate}

  Since $\,D\,$ is Noetherian, then $\,\star\,$ is  a (semi)star
operation of finite
  type on $D$.  Henceforth, it is clear  that
$\,{\calM}(\star_{_{\!f}}) = \{M\}\,.$
  \
  Thus,  $\,\widetilde{\star}$ coincides with the
  identity semistar operation $\,d_{D}\,,$\; i.e.: $$ d_D =
\widetilde{\star} \lneq \star_{_{\!f}} = \star\,. $$ Moreover, we
have already
  proved in \cite[Example 5.3]{FL3} that
  $\,\widetilde{(\star_a)} = \star_a = t_{D}\,.$  Therefore, if we
  denote by $Z$ a (new) indeterminate over the field of quotients
  $K$ of $D$, by $\calW$ the set of all rank one discrete
  valuation overrings of $D$,  then $t_D = \wedge_{\calW}$, thus:
  $$
  \begin{array}{rl}
    \Na(D, \star) =& \hskip -8pt  \Na(D, \widetilde{\star}) =\Na(D,
d_{D}) =
  D(Z)\,,\\
  \Kr(D, \star) =& \hskip -8pt \Kr(D, \star_a) = \Kr(D, t_{D})=
  \Kr(D, \wedge_{\calW})=
  \bigcap \{W(Z) \mid W \in \calW\}\,.
  \end{array}
  $$
  %%%%
    \hskip 0.5cm
 %CLAIM 1
  \bf Claim 1. \sl Let $J \in {\boldsymbol{f}}(D)$. If
$J$ is
  $\star$--e.a.b. then $J$ is principal. \rm

  Without loss of generality we can assume that $J \subseteq D$. If
  $J$ is not principal then either $J$ is not contained in any
  principal ideal of $D$ or $J$ is contained in a
  principal ideal of $D$ in both cases $J = \delta I$, for some
  nonzero ideal $I$  of $D$ such that $I^\star =   M$   and for 
some
  nonzero element $\delta \in D$ (eventually $\delta =1$). Therefore,
  $J^\star = (\delta I)^\star = \delta I^\star = \delta  M$.   
On the
  other hand, since by definition of $\star$ we have that
 $(M^{2})^\star =M^\star = M$,   then:
  $$
  (JI)^\star =  (\delta I^{2})^\star = \delta(I^\star I^\star)^\star
  =\delta (M^{2})^\star   =\delta M   = J^\star =
  (JD)^\star\,.$$
  If $J$ is
  $\star$--e.a.b. (since  $I$    is finitely generated) then 
$I^\star =
  D^\star = D$, which is a contradiction.

 %CLAIM 2
   \bf Claim 2. \sl  Let $\calV(D, \star)$ be the set of all
   $\star$--valuation overrings of $D$. If $L'$ is a
strong--$\star$--monolocality of $D$ then $L'$ is
   a localization of $D$ at a height-one prime ideal of $D$, thus
   $\calL'(D, \star) =
    \calW = \calV(D, \star)$.  In particular,
   $$\SKN(D, \star) = \bigcap \{W(Z) \mid W \in \calW\} \ (=
   \Kr(D, \star)) \,.$$
   \rm

   As a matter of fact $M = (M^{2})^\star \subseteq M^{2}L'$ and,
   since $M$ is finitely generated,  by Nakayama's Lemma we have $ML'
=
   L'$. Therefore, for some element $f \in M$, $fL'=L'$ and thus
$D_{f}
   \subseteq L'$. Since $D_{f}$ is a one-dimensional Krull domain and
   $L'$ is a  quasilocal   overring of $D_{f}$, then necessarily 
$L'$ is
   a localization of $D$ at a height-one prime ideal of $D$, hence
   $L'\in \calW$, i.e. $\calL' \subseteq \calW$.
   Conversely, if $Q$ is an height-one prime ideal of $D$, then from
   the definition of $\star$ it follows immediately that the discrete
   valuation overring
   $W:= D_{Q }$ is a strong--$\star$--monolocality of $D$. Note that
in general
   $\calL'(D, \star) \supseteq \calV(D, \star)$ and, in this case,
   each $W\in \calW$ is a $\star$--valuation overring of $D$   (by 
\cite[Theorem 3.5]{FL2}, since $W(Z) \supseteq \Kr(D, \star)$),   
that is
$\calV(D, \star) \supseteq \calW$.

 %CLAIM 3
   \bf Claim 3. \sl We have that $b_{D} \lneq \star_a
   \ (=\widetilde{(\star_a)} = t_{D})$ and $\calL'(D, b_{D} )
   \supsetneq \calL'(D,  \star_a)$ thus,  in particular,
   $$(\Kr(D, \star) = \Kr(D, \star_{a}) =)\ \SKN(D, \star_a)
   \supsetneq\SKN(D,
   b_{D}) \ (=  \Kr(D,
   b_{D}) = \Kr(D, d_{D}))\,.$$
   (cf. also Proposition
   \ref{pr:1.12} (2)).
   \rm

   We have observed above that $\star_a
   =\widetilde{(\star_a)} = t_{D}= \wedge_{\calW}$ and thus it is
   easy to see that $
   \calL'(D,  \star_a) = \calL(D,  \star_a)= \calW = \{D_{Q}\mid Q
\mbox{ is an height one
   prime ideal of } D\}$. On the other hand, $b_{D}$ coincides by
definition with
   $\wedge_{\calV}$, where  (in the present Example) we set    
$\calV :=\calV(D, d_{D}) \ ( = \calV(D,
   b_{D}))$ is the set of all the valuation
   overrings of $D$ and thus it is easily seen that  $
    \calL'(D,   \wedge_{\calV})   = \calL(D,   
\wedge_{\calV})    = \calV$.
Since there are plenty of two
   dimensional valuation overrings of $D$, then clearly $ \calL'(D,
   b_{D}) \supsetneq \calL'(D,  \star_a)$.  Finally, by
   \cite[Proposition 4.1 (5)]{FL3}, it is clear that $b_{D} \lneq
   \star_a $ if and only if     $ \Kr(D,
   b_{D})  \subsetneq   \Kr(D, \star_{a}) $.  

   %CLAIM 4
   \bf Claim 4. \sl We have that $\star \lneq \star_a
 $ (more precisely,  every nonprincipal
   ideal of $D$ is
   $\star_a$--e.a.b.,
   but not $\star$--e.a.b.) and $\calL'(D, \star )
   =\calL'(D,  \star_a)$. In particular, $\SKN(D, \star) = \SKN(D,
   \star_a)$. \rm

   Note that every
   nonzero ideal of $D$ is clearly $\star_a$--e.a.b. but from Claim 1
   we know that if an ideal of $D$ is $\star$--e.a.b. then it is a
   principal ideal. We have observed in the previous Claim 2 and 3
that
   $\calL'(D, \star ) =\calW =\calL'(D,  \star_a)$.

   %CLAIM 5
  \bf Claim 5. \sl   ${\calL}(D, \star) \supsetneq {\calL'}(D,
  \star)$ and $\KN(D, \star)  \subsetneq \SKN(D, \star)$. More
  precisely:
  $$
  \begin{array}{rl} \KN(D, \star)
  = &\hskip -8pt D(Z) \ (=\Na(D, \star)) \subsetneq \\
  \subsetneq  &\hskip -8pt  \SKN(D, \star)  =\bigcap
  \{D_{Q}(Z) \mid Q \in \Spec(D) \mbox{ and \rm{ht}} (Q)=1\} \ (=
\Kr(D, \star)) \,.
   \end{array}$$
  \rm

  By the previous considerations it is sufficient to show that
$\KN(D, \star)
  = D(Z)$.  This is an easy consequence of Claim 1, since each
  $J \in {\boldsymbol{f}}(D)$ which is a
  $\star$--e.a.b.  is a principal fractional ideal of $D$ and, by
  definition $D$ is  quasilocal   and
  $D = D^\star$, thus $D$ is a $\star$--monolocality of $D$.
  Therefore $ {\calL}(D, \star) = \{L \mid L \mbox{ is a  
quasilocal  
  overring of } $  $ D \mbox{ such }$ $\mbox{that } L =  L^{\star_{f}}
  \}$ and $D \in {\calL}(D, \star)  \setminus {\calL'}(D, \star)$ .

 %CLAIM 6
  \bf Claim 6. \sl  ${\calL}(D, \star_{a}) ={\calL'}(D,
  \star_{a}) \ (= {\calL'}(D,
  \star) = \calW)$  thus $(\KN(D, t_{D})= \ \KN(D, \star_{a})$ $ =
\SKN(D,
  \star_{a}) \ (= \SKN(D, \star) )$ and so  $\KN(D, \star) \subsetneq
  \KN(D, \star_{a})$. \rm

  This is a consequence of Corollary   \ref{cor:6.7}   since 
in this
  case we know that $\star_{a} =\widetilde{(\star_{a})}$ is a stable
(semi)star
  operation on D.

  In conluding, in this example, we have:
  $$
  d_{D}= \widetilde\star ={\star_{_{\!{\ell}}}}   ={\wedge_{\calL}}
             \lneq \star_{_{\!f}}  \lneq
             {\star_{_{\!{\ell'}}}}
={\wedge_{\calL'}}=\star_{a}=t_{D}\,.$$

  \end{exxe}

We end with an easy example announced in Remark \ref{rk:6.5}.
\medskip

%EXAMPLE 7.8
\begin{exxe} \label{ex:7.5} \sl An example of an Noetherian 
integrally closed domain $D$ with a (semi)star operation $\star$ such 
that ${\widetilde{\star}}   \lneq
           \star_{_{\!{\ell}}}$. \rm

  Let $D$ be as in Example \ref{ex:1.26}. Note that, in this case,  
$t_D = \wedge_{\calW}$ and $b_D = \wedge_{\calV}$, where $\calW $ is 
the set of all the rank 1 valuation overrings of the Krull domain $D$ 
\cite[Proposition 44.13]{G} and $\calV$ is the set of all the 
valuation overrings of $D$.  Therefore $t_D$ and $b_D$ are both 
(e.)a.b. (semi)star operations on $D$. 
  
  Let $\star := b_D$.  Since it is easy to see that $\calM(b_D) = 
\Max(D)$   \cite[Theorem 4.3 (3)]{FL3},
then $\widetilde{\ b_D} = d_D$. Moreover,  we have already seen  (in 
the proof of Claim 3, Example \ref{ex:1.26}) that
${\calL}(D, b_D)=  {\calL'}(D, b_D) =\calV$, thus 
$(b_D)_{_{\!{\ell}}} =
(b_D)_{_{\!{\ell'}}} =b_D$. 

\end{exxe}

  %%%%%%%%% BIBLIOGRAPHY

  %%%%%%%%

\end{document}